\documentclass[12pt]{amsart}

\usepackage{geometry}
\usepackage{amsfonts, adjustbox, amssymb, amscd}
\numberwithin{equation}{section}

\usepackage{bm}
\usepackage{verbatim}
\usepackage{mathrsfs}
\usepackage{graphicx}
\usepackage{tikz-cd}
\usepackage{subcaption}
\usepackage{listings}
\usepackage{subfiles}
\usepackage[toc,page]{appendix}
\usepackage{mathtools}
\usepackage{comment}
\usepackage{enumerate}
\usepackage{enumitem}
\usepackage{graphicx}
\usepackage{appendix}
\usepackage{hyperref}
\hypersetup{
    colorlinks=true,
    citecolor=red,
    linkcolor=blue,
    filecolor=magenta,      
    urlcolor=red,
}
\newcommand{\bb}{\bm{b}}
\newcommand{\Mm}{{\bf{M}}}
\newcommand{\Nn}{{\bf{N}}}

\newcommand{\Pp}{{\bf{P}}}

\newcommand{\Qq}{\mathbb{Q}}

\newcommand{\Rr}{\mathbb{R}}

\newcommand{\Center}{\operatorname{center}}

\newcommand{\Exc}{\operatorname{Exc}}

\newcommand{\ninv}{\operatorname{ninv}}
\newcommand{\inv}{\operatorname{inv}}

\newcommand{\tmld}{{\operatorname{tmld}}}

\newcommand{\Supp}{\operatorname{Supp}}

\newcommand{\mult}{\operatorname{mult}}

\newcommand{\Aa}{{\mathfrak{A}}}

\newcommand{\Ff}{\mathcal{F}}

\newcommand{\Ii}{\Gamma}

\newcommand{\Ll}{{\bf{L}}}


\newcounter{parentnumber}

\newtheorem{thm}{Theorem}[section]

\newtheorem{lem}[thm]{Lemma}
\newtheorem{prop}[thm]{Proposition}

\newtheorem{claim}[thm]{Claim}

\theoremstyle{definition}
\newtheorem{defn}[thm]{Definition}
\newtheorem{ques}[thm]{Question}
\theoremstyle{definition}
\newtheorem{rem}[thm]{Remark}

\newtheorem{deflem}[thm]{Definition-Lemma}

\newtheorem{nota}[thm]{Notation}

\theoremstyle{definition}

\begin{document}

\title{Variation of algebraically integrable adjoint foliated structures}
\author{Paolo Cascini, Jihao Liu, Fanjun Meng, Roberto Svaldi, and Lingyao Xie}

\date{\today}

\begin{abstract}
Given a canonical algebraically integrable foliation on a klt projective variety, we study the variation of the ample models of the associated adjoint foliated structures with respect to the parameter. When the foliation is of general type, we show the finiteness of ample models if the parameter is sufficiently close to $1$. When the ambient variety is of general type, we show the finiteness of ample models for all parameters. A key ingredient in our proof is the equivalence between the existence of minimal models and the termination of MMP with scaling for algebraically integrable adjoint foliated structures.
\end{abstract}

\address{Department of Mathematics, Imperial College London, 180 Queen’s
Gate, London SW7 2AZ, UK}
\email{p.cascini@imperial.ac.uk}

\address{Department of Mathematics, Peking University, No. 5 Yiheyuan Road, Haidian District, Peking 100871, China}
\email{liujihao@math.pku.edu.cn}

\address{Department of Mathematics, University of California San Diego, 9500 Gilman Drive \# 0112, La Jolla, CA 92093-0112, USA}
\email{f2meng@ucsd.edu}

\address{Dipartimento di Matematica ``F. Enriques'', Universit\`a degli Studi di Milano, Via Saldini 50, 20133 Milano (MI), Italy}
\email{roberto.svaldi@unimi.it}

\address{Department of Mathematics, University of California San Diego, 9500 Gilman Drive \# 0112, La Jolla, CA 92093-0112, USA}
\email{l6xie@ucsd.edu}

\maketitle

\pagestyle{myheadings}\markboth{\hfill P. Cascini, J. Liu, F. Meng, R. Svaldi, and L. Xie \hfill}{\hfill Variation of algebraically integrable adjoint foliated structures\hfill}

\tableofcontents

\section{Introduction}\label{sec: introduction}

We work over the field of complex numbers $\mathbb{C}$.

The existence of good minimal models for klt algebraically integrable adjoint foliated structures $(X,\Ff,t)$ of general type was established in \cite{Cas+25}. Equivalently, for $t\in[0,1)$ these are the good minimal models of the $\mathbb{R}$-divisors $tK_{\Ff}+(1-t)K_X$. Several people (e.g. D. Bejleri and K. DeVleming) asked us whether there is a wall-crossing phenomenon for such adjoint foliated structures as $t$ varies in $[0,1)$, in analogy with the wall crossing for moduli of stable pairs developed in \cite{ABIP23,MZ23}.

Since a moduli theory for adjoint foliated structures of general type has not yet been developed and several necessary tools are still unavailable, our first step toward this question is to prove the following finiteness of ample models as $t$ varies and to describe their interrelations. This will eventually yield a wall-and-chamber decomposition in the anticipated wall-crossing phenomenon. A complete wall-crossing theory (e.g. construction of reduction morphisms and wall-crossing morphisms) for adjoint foliated structures of general type requires substantial further developments and lies beyond the scope of this work. In the setting of foliations, the question becomes significantly subtler as $t\to 1$, since $K_{\Ff}$ may fail to admit a good minimal model owing to the failure of abundance for algebraically integrable foliations. Indeed, the question is closely tied to the existence of minimal models for $K_{\Ff}$, a fundamental problem in the minimal model program for foliations.

For clarity of exposition, in the following, we assume that $X$ has klt singularities, that $\Ff$ has canonical singularities, and that either $X$ or $\Ff$ is of general type, or more generally, $\lambda K_{\Ff}+(1-\lambda)K_X$ is big for some $\lambda\in [0,1]$. We provide a more general statement in Theorem \ref{thm: technical finiteness of ample models}.

\begin{thm}\label{thm: main canonical F finiteness of ample models f general type}
Let $X$ be a klt projective variety and $\Ff$ a canonical algebraically integrable foliation on $X$. For $t\in[0,1]$ set $\Aa_t:=(X,\Ff,t)$. Assume that $\Aa_{\lambda}$ is of general type for some $\lambda\in [0,1]$. Then there exist a rational number $\epsilon<\lambda$ and finitely many rational numbers
$$\max\{\epsilon,0\}:=t_1<t_2<\dots<t_n:=1$$
with the following properties. Let $\Ii:=\{\{t_i\},\,(t_i,t_{i+1})\mid 1\leq i\leq n-1\}$. For any $\mathcal{P}\in\Ii$ there exists a rational map $\psi_{\mathcal{P}}\colon X\dashrightarrow Z_{\mathcal{P}}$ such that:
\begin{enumerate}
    \item For any $\mathcal{P}\in\Ii$ and $t\in\mathcal{P}$, the map $\psi_{\mathcal{P}}$ is the ample model of $\Aa_t$.
    \item For any $\mathcal{P},\mathcal{Q}\in\Ii$ with $\mathcal{Q}\subset\partial\mathcal{P}$, there exists a unique contraction $\mu_{\mathcal{P},\mathcal{Q}}\colon Z_{\mathcal{P}}\to Z_{\mathcal{Q}}$ such that $\mu_{\mathcal{P},\mathcal{Q}}\circ\psi_{\mathcal{P}}=\psi_{\mathcal{Q}}$.
\end{enumerate}
\end{thm}

Compared with the case of pairs, a key difference is that $K_{\Ff}$ may not admit a good minimal model (cf. \cite[Theorem 3 IV.5.11]{McQ08}, \cite[Subsection 5.4]{ACSS21}). Consequently, we work over the half-open interval $[\epsilon,1)$, whereas classical results for pairs typically hold over a closed interval.

A central ingredient in the proof of Theorem \ref{thm: main canonical F finiteness of ample models f general type} is an equivalence between the existence of minimal models and the termination of the MMP with scaling for algebraically integrable adjoint foliated structures. 

\begin{thm}\label{thm: scaling vs mmp main}
Let $(X,\Ff,t)$ be a projective lc algebraically integrable adjoint foliated structure with $X$ klt. If $(X,\Ff,t)$ admits a minimal model, then one can run a $(tK_{\Ff}+(1-t)K_X)$-MMP with scaling of an ample divisor on $X$, and any such MMP terminates with a minimal model of $(X,\Ff,t)$.
\end{thm}

We also expect Theorem \ref{thm: scaling vs mmp main} to be useful for further developments in the minimal model program for algebraically integrable adjoint foliated structures. More general versions of Theorem \ref{thm: scaling vs mmp main} are given in Theorems \ref{thm: eomm implies tof with scaling} and \ref{thm: eo bswlcm implies eomm}.

It is interesting to ask that, in the setting of Theorem \ref{thm: main canonical F finiteness of ample models f general type}, whether we can choose $\epsilon$ and $\{t_i\}$ to be independent of $X$ or $\Ff$. More precisely, we ask the following questions:
\begin{ques}
    Notation and conditions as in \ref{thm: main canonical F finiteness of ample models f general type}.
    \begin{enumerate}
        \item Assume that $\dim X=d$ is fixed. 
        \begin{enumerate}
            \item Can we choose $\epsilon=\epsilon(d,\lambda)$ so that it only depends on $d$ and $\lambda$?
            \item Can we find an explicit or the optimal value of $\epsilon(d,1)$?
        \end{enumerate}
        \item Assume that $X$ is fixed and $\lambda=0$, i.e. $X$ is of general type. We have $t_1=0$. 
        \begin{enumerate}
            \item Can we choose $t_2$ so that it only depends on $X$?
            \item Can we choose $n$ so that it only depends on $X$? What about the set $\{t_i\}_{i=1}^{n}$?
        \end{enumerate}
    \end{enumerate}
\end{ques}

\medskip

\noindent\textit{Structure of the paper.} In Section \ref{sec: preliminaries} we recall some preliminary results and prove some perturbation results that will be used later. In Section \ref{sec: models} we define different types of models, particularly models in the sense of Birkar-Shokurov, of adjoint foliated structures, and study their basic behaviors. In Section \ref{sec: mmp scaling} we show the equivalence between existence of minimal models and termination of MMP with scaling, and prove Theorem \ref{thm: scaling vs mmp main}. In Section \ref{sec: Finiteness of ample models} we prove a finiteness of ample models result which implies Theorem \ref{thm: main canonical F finiteness of ample models f general type}.

\medskip 

\noindent\textbf{Acknowledgements.} The authors would like to thank Dori Bejleri and Kristin DeVleming for asking the question. They would like to thank Calum Spicer and Nikolaos Tsakanikas for useful discussions. A part of the work was done when the authors attended the 2025 SRI conference and they would like to thank the conference organizers for their hospitality. 

The work is supported by the National Key R\&D Program of China \#2024YFA1014400. The first author is partially funded by a Simons Collaboration Grant. The third author is partially supported by a Simons Collaboration Grant and an AMS-Simons Travel Grant. The fourth author is supported by the “Programma per giovani ricercatori Rita Levi Montalcini” of the Italian Ministry of University and Research and by PSR 2022 – Linea 4 of the University of Milan. He is a member of the GNSAGA group of INDAM. The fifth author is supported by a grant from the Simons Foundation.

\section{Preliminaries}\label{sec: preliminaries}

We will adopt the standard notation and definitions on MMP in \cite{KM98,BCHM10} and use them freely. For adjoint foliated structures, generalized foliated quadruples, foliated triples, foliations, and generalized pairs, we adopt the notation and definitions in \cite{Cas+24,Cas+25} (for adjoint foliated structures) which generally align with \cite{LLM23,CHLX23} (for generalized foliated quadruples), \cite{CS20, ACSS21, CS21} (for foliations and foliated triples), and \cite{BZ16,HL23} (for generalized pairs and $\bb$-divisors), possibly with minor differences. 

\subsection{Notation}

\begin{defn}
    A \emph{contraction} is a projective morphism of varieties
    $f \colon X \to Y$ such that $f_\ast \mathcal{O}_X=\mathcal{O}_Y$.
\end{defn}

\begin{nota}
    Let $h \colon  X\dashrightarrow X'$ be a birational map between normal varieties. We denote by $\Exc(h)$ the reduced divisor supported on the codimension one part of the exceptional locus of $h$. 
\end{nota}

\begin{defn}
    Let $X\rightarrow U$ be a projective morphism from a normal quasi-projective variety to a variety.  Let $D$ be an $\Rr$-Cartier $\Rr$-divisor on $X$ and $\phi: X\dashrightarrow X'$ a birational map$/U$. We say that $\phi$ is $D$-negative (resp. $D$-trivial) if the following conditions hold:
    \begin{enumerate}
    \item $\phi$ does not extract any divisor.
    \item $D':=\phi_\ast D$ is $\Rr$-Cartier.
    \item There exists a resolution of indeterminacy $p: W\rightarrow X$ and $q: W\rightarrow X'$, such that
    $$p^\ast D=q^\ast D'+F$$
    where $F\geq 0$ and $\Supp p_\ast F=\Exc(\phi)$ (resp. $F=0$).
    \end{enumerate}
\end{defn}

\begin{defn}[NQC]
    Let $X\rightarrow U$ be a projective morphism between normal quasi-projective varieties. Let $D$ be a nef $\Rr$-Cartier $\Rr$-divisor on $X$ and $\Mm$ a nef $\bb$-divisor on $X$. We say that $D$ is \emph{NQC}$/U$ if $D=\sum d_iD_i$, where each $d_i\geq 0$ and each $D_i$ is a nef$/U$ Cartier divisor. We say that $\Mm$ is \emph{NQC}$/U$ if $\Mm=\sum \mu_i\Mm_i$, where each $\mu_i\geq 0$ and each $\Mm_i$ is a nef$/U$ $\bb$-Cartier $\bb$-divisor.
\end{defn}

\begin{defn}
    Let $\pi\colon X\rightarrow U$ be a projective morphism from a normal quasi-projective variety to a quasi-projective variety, $D$ a pseudo-effective$/U$ $\Rr$-Cartier $\Rr$-divisor on $X$, and $P$ a prime divisor on $X$. We define $\sigma_{P}(X/U,D)$ as in \cite[Definition 3.1]{LX25} by considering $\sigma_{P}(X/U,D)$ as a number in  $[0,+\infty)\cup\{+\infty\}$. We define $N_{\sigma}(X/U,D)=\sum_Q\sigma_Q(X/U,D)Q$
    where the sum runs through all prime divisors on $X$ and consider it as a formal sum of divisors with coefficients in $[0,+\infty)\cup\{+\infty\}$. We say that $D$ is \emph{movable$/U$} if $N_{\sigma}(X/U,D)=0$.
\end{defn}

\subsection{Adjoint foliated structures}

\begin{defn}[Foliations, {cf. \cite{ACSS21,CS21}}]\label{defn: foliation}
Let $X$ be a normal variety. A \emph{foliation} on $X$ is a coherent subsheaf $T_{\Ff}\subset T_X$ such that
\begin{enumerate}
    \item $T_{\Ff}$ is saturated in $T_X$, i.e. $T_X/T_{\Ff}$ is torsion free, and
    \item $T_{\Ff}$ is closed under the Lie bracket.
\end{enumerate}
The \emph{canonical divisor} of $\Ff$ is a divisor $K_\Ff$ such that $\mathcal{O}_X(-K_{\mathcal{F}})\cong\mathrm{det}(T_\Ff)$. If $T_\Ff=0$, then we say that $\Ff$ is a \emph{foliation by points}.

Given any dominant map 
$h: Y\dashrightarrow X$ and a foliation $\mathcal F$ on $X$, we denote by $h^{-1}\Ff$ the \emph{pullback} of $\Ff$ on $Y$ as constructed in \cite[3.2]{Dru21}. Given any birational map $g: X\dashrightarrow X'$, we denote by $g_\ast \Ff:=(g^{-1})^{-1}\Ff$ the \emph{pushforward} of $\Ff$ on $X'$. We say that $\Ff$ is an \emph{algebraically integrable foliation} if there exists a dominant map $f: X\dashrightarrow Z$ such that $\Ff=f^{-1}\Ff_Z$, where $\Ff_Z$ is the foliation by points on $Z$, and we say that $\Ff$ is \emph{induced by} $f$.

A subvariety $S\subset X$ is called \emph{$\Ff$-invariant} if for any open subset $U\subset X$ and any section $\partial\in H^0(U,\Ff)$, we have $\partial(\mathcal{I}_{S\cap U})\subset \mathcal{I}_{S\cap U}$,  where $\mathcal{I}_{S\cap U}$ denotes the ideal sheaf of $S\cap U$ in $U$.  
For any prime divisor $P$ on $X$, we define $\epsilon_{\Ff}(P):=1$ if $P$ is not $\Ff$-invariant and $\epsilon_{\Ff}(P):=0$ if $P$ is $\Ff$-invariant. For any prime divisor $E$ over $X$, we define $\epsilon_{\Ff}(E):=\epsilon_{\Ff_Y}(E)$ where $h: Y\dashrightarrow X$ is a birational map such that $E$ is on $Y$ and $\Ff_Y:=h^{-1}\Ff$.

Suppose that the foliation structure $\Ff$ on $X$ is clear in the context. Then, given an $\mathbb R$-divisor $D =\sum a_iD_i$ where each $D_i$ is a prime divisor,
we denote by $D^{\ninv} \coloneqq \sum \epsilon_{\Ff}(D_i)a_iD_i$ and $D^{\inv} \coloneqq D-D^{\ninv}$.
\end{defn}

\begin{defn}\label{defn: afs}
An \emph{adjoint foliated structure} $\Aa/U:=(X,\Ff,B,\Mm,t)/U$ is the datum of a normal quasi-projective variety $X$ and a projective morphism $X\rightarrow U$, a foliation $\Ff$ on $X$, an $\Rr$-divisor $B\geq 0$ on $X$, a $\bb$-divisor $\Mm$ nef$/U$, and a real number $t\in [0,1]$ such that $K_{\Aa}:=tK_{\Ff}+(1-t)K_X+B+\Mm_X$ is $\Rr$-Cartier. We may simply say that ``$\Aa/U$ is an adjoint foliated structure" without mentioning $X,\Ff,B,\Mm,t$ at all. $X,t$ are called the \emph{ambient variety} and \emph{parameter} of $\Aa$ respectively. 

We say that $\Aa/U$ is \emph{of general type} if $K_{\Aa}$ is big$/U$. We say that $\Aa$ is \emph{algebraically integrable} if $\Ff$ is algebraically integrable. We say that $\Aa/U$ is \emph{KNQC} if $K_{\Aa}$ is NQC$/U$.

For any $\Rr$-divisor $D$ on $X$ and nef$/U$ $\bb$-divisor $\Nn$ on $X$ such that $D+\Nn_X$ is $\Rr$-Cartier, we denote by $(\Aa,D,\Nn):=(X,\Ff,B+D,\Mm+\Nn,t)$. If $D=0$ then we may drop $D$, and if $\Nn=\bm{0}$ then we may drop $\Nn$.

When $t=0$ or $\Ff=T_X$, we call $(X,B,\Mm)/U$ a \emph{generalized pair}, and in addition, if $\Mm=\bm{0}$, then we call $(X,B)/U$ a pair. If $B=0$, or if $\Mm=\bm{0}$, or if $U$ is not important, then we may drop $B,\Mm,U$ respectively. If $U=\{pt\}$ then we also drop $U$ and say that $(X,\Ff,B,\Mm,t)$ is \emph{projective}.

For any birational map$/U$ $\phi: X\dashrightarrow X'$, we define $\phi_*\Aa:=(X',\phi_*\Ff,\phi_*B,\Mm,t)$ and say that $\phi_*\Aa$ is the \emph{image} of $\Aa$ on $X'$. 
For any projective birational morphism $h: X'\rightarrow X$, we define 
$$h^*\Aa:=(X',\Ff',B',\Mm,t)$$
where $\Ff':=h^{-1}\Ff$ and $B'$ is the unique $\Rr$-divisor such that $K_{h^*\Aa}=h^*K_{\Aa}$. For any prime divisor $E$ on $X'$, we denote by
$$a(E,\Aa):=-\mult_EB'$$
the \emph{discrepancy} of $E$ with respect to $\Aa$. The \emph{total minimal log discrepancy} of $\Aa$ is
$$\tmld(\Aa):=\inf\{a(E,\Aa)+t\epsilon_{\Ff}(E)+(1-t)\mid E\text{ is over }X\}.$$
For any non-negative real number $\epsilon$, we say that $\Aa$ is \emph{$\epsilon$-lc} (resp. \emph{$\epsilon$-klt}) if $\tmld(\Aa)\geq\epsilon$ (resp. $>\epsilon$). We say that $\Aa$ is \emph{lc} (resp. \emph{klt}) if $\Aa$ is $0$-lc (resp. $0$-klt). An \emph{lc place} of $\Aa$ is a prime divisor $E$ over $X$ such that $a(E,\Aa)=-t\epsilon_{\Ff}(E)-(1-t)$. An \emph{lc center} of $\Aa$ is the image of an lc place of $\Aa$ on $X$.

For any adjoint foliated structures $\Aa_i/U=(X,\Ff,B_i,\Mm_i,t_i)/U$ and real numbers $a_i\in [0,1]$ such that $\sum a_i=1$, we denote by
$$\sum a_i\Aa_i:=\left(X,\Ff,\sum a_iB_i,\sum a_i\Mm_i,\sum a_it_i\right).$$
\end{defn}

\begin{defn}[Qdlt, {\cite[Definition 3.1]{Cas+24}}]\label{defn: qdlt afs}
    We say that an adjoint foliated structure $\Aa/U:=(X,\Ff,B^{\ninv}+(1-t)B^{\inv},\Mm,t)/U$ is \emph{qdlt} if $t<1$, $\Aa$ is lc, $(X,B,\Mm)$ is qdlt (cf. \cite[Definition 7.1.1]{CHLX23}), and any lc place of $\Aa$ is an lc place of $(X,B,\Mm)$.
\end{defn}

\begin{deflem}\label{deflem: eoqdlt model}
    Let $\Aa/U$ be an lc algebraically integrable adjoint foliated structure with ambient variety $X$ and parameter $t<1$. A \emph{$\mathbb Q$-factorial qdlt modification} of $\Aa$ is a birational morphism $h: X'\rightarrow X$ such that $h^*\Aa$ is $\mathbb Q$-factorial qdlt and $h$ only extracts lc places of $\Aa$, which always exists by \cite[Theorem 1.9]{Cas+24}.
\end{deflem}

\begin{defn}[Potentially klt]\label{defn: potentially klt}
Let $X$ be a normal quasi-projective variety. We say that $X$ is \emph{potentially klt} if $(X,\Delta)$ is klt for some $\Rr$-divisor $\Delta\geq 0$. 
\end{defn}

\begin{defn}[{\cite[Definition 3.11]{Cas+25}}]\label{defn: foliated log smooth}
Let $\Aa/U:=(X,\Ff,B,\Mm,t)/U$ be an algebraically integrable adjoint foliated structure. We say that $\Aa$ is \emph{foliated log smooth} if there exists a contraction $f: X\rightarrow Z$ satisfying the following.
\begin{enumerate}
  \item $X$ has at most quotient singularities.
  \item $\Ff$ is induced by $f$.
  \item $(X,\Sigma_X)$ is toroidal for some reduced divisor $\Sigma_X$ such that $\Supp B\subset\Sigma_X$.  In particular, $(X,\Supp B)$ is toroidal, and $X$ is $\Qq$-factorial klt.
  \item There exists a log smooth pair $(Z,\Sigma_Z)$ such that $$f: (X,\Sigma_X)\rightarrow (Z,\Sigma_Z)$$ is an equidimensional toroidal contraction.
  \item $\Mm$ descends to $X$.
\end{enumerate}
We say that $f: (X,\Sigma_X)\rightarrow (Z,\Sigma_Z)$ is \emph{associated with} $(X,\Ff,B,\Mm,t)$.
\end{defn}

\begin{defn}[Foliated log resolutions]\label{defn: log resolution}
Let $\Aa/U:=(X,\Ff,B,\Mm,t)/U$ be an algebraically integrable adjoint foliated structure. A \emph{foliated log resolution} of $\Aa$ is a birational morphism $h: X'\rightarrow X$ such that 
$$(X',\Ff':=h^{-1}\Ff,B':=h^{-1}_\ast B+\Exc(h),\Mm,t)$$ 
is foliated log smooth. By \cite[Lemma 6.2.4]{CHLX23}, a foliated log resolution for $\Aa$ always exists.
\end{defn}

\begin{defn}[Property $(\ast )$ foliations, {\cite[Definition 3.8]{ACSS21}, \cite[Definition 7.2.2]{CHLX23}}]\label{defn: foliation property *}
Let $(X,\Ff,B,\Mm)/U$ be a generalized foliated quadruple, $G\geq 0$ be a reduced divisor on $X$, and $f: X\rightarrow Z$ a contraction. We say that $(X,\Ff,B,\Mm;G)/Z$ satisfies \emph{Property $(\ast )$} if the following conditions hold.
\begin{enumerate}
\item $\Ff$ is induced by $f$ and $G$ is an $\Ff$-invariant divisor.
\item $f(G)$ is of pure codimension $1$, $(Z,f(G))$ is log smooth, and $G=f^{-1}(f(G))$.
\item For any closed point $z\in Z$ and any reduced divisor  $\Sigma\ge f(G)$ on $Z$ such that  $(Z,\Sigma)$ is log smooth near $z$, $(X,B+G+f^\ast (\Sigma-f(G)),\Mm)$ is lc over a neighborhood of $z$.
\end{enumerate}
We say that $f$, $Z$, and $G$ are \emph{associated} with $(X,\Ff,B,\Mm)$. 
\end{defn}

\begin{defn}[ACSS, {cf. \cite[Definitions 5.4.2, 7.2.2, 7.2.3]{CHLX23}}]\label{defn: ACSS f-triple}
Let $(X,\Ff,B,\Mm)$ be an lc generalized foliated quadruple, $G\geq 0$ a reduced divisor on $X$, and $f: X\rightarrow Z$ a contraction. We say that $(X,\Ff,B,\Mm;G)/Z$ is \emph{ACSS} if the following conditions hold:
\begin{enumerate}    
\item $(X,\Ff,B,\Mm;G)/Z$ satisfies Property $(\ast )$.
\item $f$ is equidimensional.
\item There exists an $\Rr$-Cartier $\Rr$-divisor $D\geq 0$ on $X$ and a nef$/X$ $\bb$-divisor $\Nn$ on $X$, such that  $\Supp\{B\}\subset\Supp D$, $\Nn-\alpha\Mm$ is nef$/X$ for some $\alpha>1$, and for any reduced divisor $\Sigma\geq f(G)$ such that $(Z,\Sigma)$ is log smooth, $$(X,B+D+G+f^\ast (\Sigma-f(G)),\Nn)$$ 
      is qdlt (cf. \cite[Definition 7.1.1]{CHLX23}, \cite[Definition 35]{dFKX17}).
\item For any lc center of $(X,\Ff,B,\Mm)$ with generic point $\eta$, over a neighborhood of $\eta$,
    \begin{enumerate}
      \item $\eta$ is the generic point of an lc center of $(X,\Ff,\lfloor B\rfloor)$, and
       \item $f: (X,B+G)\rightarrow (Z,f(G))$ is a toroidal morphism.
    \end{enumerate}
\end{enumerate}
If $(X,\Ff,B,\Mm;G)/Z$ is ACSS, then we say that $(X,\Ff,B,\Mm)/Z$ and $(X,\Ff,B,\Mm)$ are ACSS.
\end{defn}

\begin{defn}[ACSS modifications]\label{defn: simple model}
    Let $\Aa=(X,\Ff,B,\Mm)$ and $\Aa'=(X',\Ff',B',\Mm)$ be two lc algebraically integrable generalized foliated quadruples and let $h: X'\rightarrow X$ be a birational morphism. We say that $h: X'\rightarrow X$ is a \emph{$\mathbb Q$-factorial ACSS modification} of $\Aa$ if there exists a contraction $f: X'\rightarrow Z$ and a reduced divisor $G$ on $X'$, such that
    \begin{enumerate}
        \item $\Aa'=h^*\Aa$,
        \item $(X',\Ff',B',\Mm;G)/Z$ is $\mathbb Q$-factorial ACSS, and
        \item  $a(E,\Aa)=-\epsilon_{\Ff}(E)$ for any $h$-exceptional prime divisor $E$.
    \end{enumerate}
    $\mathbb Q$-factorial ACSS modifications always exist by \cite[Theorem 8.2.2]{CHLX23}.
\end{defn}

\section{Models of adjoint foliated structures}\label{sec: models}

This section aligns with \cite[Section 2]{Bir12}, \cite[Section 3]{HL23}, and \cite[Section 4]{LMX24} and the proofs are similar.

\begin{defn}[Log birational models]\label{defn: log birational model}
Let $\Aa/U$ be an adjoint foliated structure with ambient variety $X$, $\phi: X\dashrightarrow X'$ a birational map$/U$, and $E:=\Exc(\phi^{-1})$ the reduced $\phi^{-1}$-exceptional divisor. Assume that $a(D,\Aa)\leq-t\epsilon_{\Ff}(D)-(1-t)$ for any component $D$ of $E$. We let
  $$\Aa':=\left(\phi_*\Aa,-\sum_Da(D,\Aa)D\right)$$
where the sum runs through all irreducible components of $E$. If $K_{\Aa'}$ is $\mathbb R$-Cartier then we say that $\Aa'/U$ is a \emph{log birational model} of $\Aa/U$.
\end{defn}

\begin{defn}[Minimal models]\label{defn: minimal model}
    Let $\Aa/U$ be an adjoint foliated structure with ambient variety $X$ and $\Aa'/U$ a log birational model of $\Aa/U$ with ambient variety $X'$ and associated  birational map $\phi: X\dashrightarrow X'$, such that $K_{\Aa'}$ is nef$/U$. 
    %Let $X,X'$ be the ambient varieties of $\Aa$ and $\Aa'$ respectively and 
    Let $t$ be the parameter of $\Aa$.
    \begin{enumerate}
        \item We say that $\Aa'/U$ is a \emph{bs-weak lc model} or \emph{weak lc model in the sense of Birkar-Shokurov} of $\Aa/U$, if for any prime divisor $D$ on $X$ which is exceptional over $X'$, $$a(D,\Aa)\leq a(D,\Aa').$$
        We also say that $\phi$ is a bs-weak lc model of $\Aa/U$.
        \item We say that $\Aa'/U$ is a \emph{bs-minimal model} or \emph{minimal model in the sense of Birkar-Shokurov} of $\Aa/U$, if for any prime divisor $D$ on $X$ which is exceptional over $X'$, $$a(D,\Aa)<a(D,\Aa').$$
        We also say that $\phi$ is a bs-minimal model of $\Aa/U$.
        \item We say that $\Aa'/U$ is a \emph{bs-semi-ample model} or \emph{semi-ample model in the sense of Birkar-Shokurov} of $\Aa/U$ if it is a bs-weak lc model of $\Aa/U$ and $K_{\Aa'}$ is semi-ample$/U$.      We also say that $\phi$ is a bs-semi-ample model of $\Aa/U$.
        \item We say that $\Aa'/U$ is a \emph{bs-good minimal model} or \emph{good minimal model in the sense of Birkar-Shokurov} of $\Aa/U$ if it is a bs-minimal model of $\Aa/U$ and $K_{\Aa'}$ is semi-ample$/U$. We also say that $\phi$ is a bs-good minimal model of $\Aa/U$.
        \end{enumerate}
If, in addition, the induced birational map $X\dashrightarrow X'$ does not extract any divisor, then we remove the initial ``bs-" or the phrase ``in the sense of Birkar-Shokurov" in the previous definitions. 
\begin{enumerate}
        \item[(5)] We say that $\Aa'/U$ is a \emph{log minimal model} of $\Aa/U$ if it is a bs-minimal model of $\Aa/U$, $\Aa'$ is $\Qq$-factorial qdlt if $t<1$, and $\Aa'$ is $\Qq$-factorial ACSS if $t=1$. We also say that $\phi$ is a log minimal model of $\Aa/U$.
        \item[(6)] We say that $\Aa'/U$ is a \emph{good log minimal model} of $\Aa/U$ if it is a log minimal model of $\Aa$ and $K_{\Aa'}$ is semi-ample$/U$. We also say that $\phi$ is a good log minimal model of $\Aa/U$.
\end{enumerate}
\end{defn}

\begin{lem}\label{lem: hl23 2.6}
Let $\Aa/U$ be an adjoint foliated structure and let $\Aa'/U$ be a bs-weak lc model of $\Aa/U$. Let $X,X'$ be the ambient varieties of $\Aa$ and $\Aa'$ respectively and $\phi: X\dashrightarrow X'$ the associated birational map. Let $p: W\rightarrow X$ and $q: W\rightarrow X'$ be birational morphisms such that $q=\phi\circ p$. Assume that
$$p^*K_{\Aa}=q^*K_{\Aa'}+E,$$
then $E\geq 0$ and is exceptional$/X'$.
\end{lem}
\begin{proof}
For any prime divisor $D$ that is an irreducible component of $E$, $$\mult_DE=a(D,\Aa')-a(D,\Aa).$$ 
Therefore, if $D$ is not exceptional$/X$, then:
\begin{itemize}
    \item If $D$ is not exceptional$/X'$, then $\mult_DE=0$ by Definition \ref{defn: log birational model}.
    \item If $D$ is exceptional$/X'$, then $\mult_DE\geq 0$ by Definition \ref{defn: minimal model}(1).
\end{itemize}
Therefore, $p_*E\geq 0$. Since $K_{\Aa'}$ is nef$/U$, $q^*K_{\Aa'}$ is nef$/X$, hence $E$ is anti-nef$/X$. By the negativity lemma, $E\geq 0$.

If $E$ is not exceptional$/X'$, then there exists a component $D$ of $E$ that is not exceptional$/X'$. If $D$ is not exceptional$/X$, then $\mult_DE=0$ by Definition \ref{defn: log birational model}, a contradiction. Thus $D$ is exceptional over $X$. In particular, $\phi$ extracts $D$. Since $\Aa'/U$ is a log birational model of $\Aa$, $$a(D,\Aa')=a(D,\Aa),$$ 
which implies that $\mult_DE=0$, a contradiction.
\end{proof}

\begin{lem}\label{lem: g-pair version bir12 2.7}
Let $\Aa/U$ be an adjoint foliated structure. Let $\Aa_1/U$ and $\Aa_2/U$ be two bs-weak lc models of $\Aa/U$. Let $X,X_1,X_2$ be the ambient varieties of $\Aa,\Aa_1,\Aa_2$ respectively with induced birational maps $\phi: X_1\dashrightarrow X_2$. Let $h_1: W\rightarrow X_1$ and $h_2: W\rightarrow X_2$ be two birational morphisms such that $\phi\circ h_1=h_2$. Then:
\begin{enumerate}
    \item $h_1^*K_{\Aa_1}=h_2^*K_{\Aa_2}.$
    \item If $K_{\Aa_2}$ is semi-ample$/U$, then $K_{\Aa_1}$ is semi-ample$/U$.
    \item If $K_{\Aa_2}$ is ample$/U$, then $\phi$ is a morphism.
\end{enumerate}
\end{lem}
\begin{proof}
Let $\phi_1: X\dashrightarrow X_1$ and $\phi_2: X\dashrightarrow X_2$ be the induced birational maps. Possibly replacing $W$ with a higher model, we may assume that the induced birational map $h: W\rightarrow X$ is a morphism. Let $$E_i:=h^*K_{\Aa}-h_i^*K_{\Aa_i}$$
for $i\in\{1,2\}$. By Lemma \ref{lem: hl23 2.6}, $E_i\geq 0$ and is exceptional over $X_i$ for $i\in\{1,2\}$. Thus $h_{1,*}(E_2-E_1)\geq 0$ and $E_1-E_2$ is nef$/X_1$, and $h_{2,*}(E_1-E_2)\geq 0$ and $E_2-E_1$ is nef$/X_2$. By the negativity lemma, $E_2-E_1\geq 0$ and $E_1-E_2\geq 0$. Thus $E_1=E_2$, which implies (1). (2) immediately follows from (1). By (1), if $K_{\Aa_2}$ is ample$/U$, then $h_2: W\rightarrow X_2$ is the ample model$/U$ of $h^*K_{\Aa_1}$, hence $\phi$ is the ample model$/U$ of $K_{\Aa_1}$. Since $K_{\Aa_1}$ is semi-ample$/U$, $\phi$ is a morphism. This implies (3).
\end{proof}

\begin{lem}\label{lem: numerical equivalence model}
    Let $r$ be a positive real number. Let $\Aa_1/U$ and $\Aa_2/U$ be two adjoint foliated structures with the same ambient variety $X$ and  such that
    $$K_{\Aa_2}\equiv_U rK_{\Aa_1}.$$
    Let $\Aa_1'/U$ be a weak lc model (resp. minimal model) of $\Aa_1/U$ with ambient variety $X'$ and induced birational map $\phi: X\dashrightarrow X'$. Let $\Aa_2':=\phi_*\Aa$. Then:
    \begin{enumerate}
        \item $\Aa_2'/U$ is a weak lc model (resp. minimal model) of $\Aa_2/U$.
        \item If $\Aa_1'/U$ is a semi-ample model (resp. good minimal model) of $\Aa_1/U$ and
    $$K_{\Aa_2}\sim_{\mathbb R,U} rK_{\Aa_1},$$
      then  $\Aa_2'/U$ is a semi-ample model (resp. good minimal model) of $\Aa_2/U$.
    \end{enumerate}
\end{lem}
\begin{proof}
Let $p: W\rightarrow X$ and $q: W\rightarrow X'$ be a resolution of indeterminacy. By Lemma \ref{lem: hl23 2.6},
$$p^*K_{\Aa_1}=q^*K_{\Aa_1'}+E$$
for some $\Rr$-divisor $E\geq 0$ that is exceptional$/X'$. Then
$$p^*K_{\Aa_2}\equiv_U rq^*K_{\Aa_1'}+rE,$$
so that
$$K_{\Aa_2'}=q_*p^*K_{\Aa_2}\equiv q_*(rq^*K_{\Aa_1'}+rE)=rK_{\Aa_1'}$$
is nef$/U$. Moreover, if $K_{\Aa_1'}$ is semi-ample$/U$ and $K_{\Aa_2}\sim_{\mathbb R,U} rK_{\Aa_1}$, then 
$$K_{\Aa_2'}\sim_{\mathbb R,U}rK_{\Aa_1'}$$
is semi-ample$/U$. 

We have
$$p^*K_{\Aa_2}\equiv_U q^*K_{\Aa_2'}+rE.$$
Therefore, for any prime divisor $D$ on $X$ which is exceptional over $X'$,
$$a(D,\Aa_2')-a(D,\Aa_2)=-\mult_Dp_*(rE)=r(a(D,\Aa_1')-a(D,\Aa_1)).$$
Thus, 
$a(D,\Aa_2)\leq\text{(resp. }<\text{) }a(D,\Aa_2')$ if and only if $a(D,\Aa_1)\leq\text{(resp. }<\text{) }a(D,\Aa_1')$ and the lemma follows immediately from the definitions.
\end{proof}

\begin{thm}[Very exceptional MMP]\label{thm: very exceptional mmp}
Let $\Aa/U$ be an lc algebraically integrable adjoint foliated structure with ambient variety $X$ and parameter $t$. Assume that $X$ is potentially klt, and $K_{\Aa}\sim_{\mathbb R,U}E\geq 0$ for some very exceptional$/U$ (cf. \cite[Definition 3.1]{Bir12}) $\mathbb R$-divisor $E$. Then we may run a $K_{\Aa}$-MMP$/U$ with scaling of an ample$/U$ $\Rr$-divisor and any such MMP terminates with a good minimal model $\Aa'/U$ such that $K_{\Aa'}\sim_{\mathbb R,U}0$.
\end{thm}
\begin{proof}
    Let $\phi_i: \Aa_i\dashrightarrow \Aa_{i+1}$ be a step of a $K_{\Aa}$-MMP$/U$ with scaling of an ample$/U$ $\Rr$-divisor $A$. For each $i$, let $X_i$ the ambient variety of $\Aa_i$, let $A_i$ and $E_i$ be the images of $A$ and $E$ on $X_i$ respectively, and let 
    $$\lambda_i:=\inf\{s\geq 0\mid K_{\Aa_i}+sA_i\text{ is nef}\}$$
    be the scaling numbers. By \cite[Theorem 8.1]{Cas+25}, either this MMP terminates with $\Aa_n/U$ for some $n\geq 0$, or $\lim_{i\rightarrow+\infty}\lambda_i=0$. In the latter case, we let $n$ be a positive integer such that $\phi_i$ is a flip for any $i\geq n$, hence
    $$K_{\Aa_n}=\lim_{i\rightarrow+\infty}(\psi_{i,n})_*\left(K_{\Aa_i}+\lambda_iA_i\right)$$
    is movable$/U$, where $\psi_{i,n}: X_i\dashrightarrow X_n$ is the induced birational map. In either case, $K_{\Aa_n}\sim_{\mathbb R}E_n\geq 0$ is movable$/U$ and very exceptional$/U$, hence $E_n=0$ by \cite[Lemma 3.4]{Bir12}. Thus, we have $\Aa'=\Aa_n$.
\end{proof}

\begin{defn}[Foliated log smooth model]\label{defn: log smooth models}
Let $\Aa/U:=(X,\Ff,B,\Mm,t)/U$ be an lc algebraically integrable adjoint foliated structure and $h: X'\rightarrow X$ a foliated log resolution of $\Aa$. Let $E\geq 0$ be an $h$-exceptional $\mathbb R$-divisor and let $\Aa':=(h^*\Aa,E)$. We say that $\Aa'$ is a \emph{foliated log smooth model} of $\Aa$ if
\begin{enumerate}
    \item $\Aa'$ is foliated log smooth and lc, and
    \item for any $h$-exceptional prime divisor $D$ that is not an lc place of $\Aa$, $D$ is a component of $E$.
\end{enumerate}
\end{defn}

\begin{lem}\label{lem: g-pair version bir12 2.8}
Let $\Aa/U$ be an lc algebraically integrable adjoint foliated structure with ambient variety $X$. Let $\Aa_W$ be a foliated log smooth model of $\Aa$ with ambient variety $W$. 

Then any bs-weak lc model (resp. bs-minimal model, bs-semi-ample model, bs-good minimal model, log minimal model, good log minimal model) of $\Aa_W/U$ is a bs-weak lc model (resp. bs-minimal model, bs-semi-ample model, bs-good minimal model, log minimal model, good log minimal model) of $\Aa/U$. 
\end{lem}

\begin{proof}
We let $h: W\rightarrow X$ be the induced birational morphism. We may write
$$K_{\Aa_W}=h^*K_{\Aa}+E$$
for some $E\geq 0$ that is $h$-exceptional, and $D\subset\Supp E$ for any $h$-exceptional prime divisor $D$ that is not an lc place of $\Aa$.

\begin{claim}\label{claim: log smooth model log discrepancy compare}
Let $\Aa'/U$ be a bs-weak lc model of $\Aa_W/U$. Then $$a(D,\Aa)\leq a(D,\Aa')$$ for any prime divisor $D$ over $X$.
\end{claim}
\begin{proof}
Let $X'$ be the ambient variety of $\Aa'$, $\phi_W: W\dashrightarrow X'$ be the induced birational map, and let $p: V\rightarrow W$ and $q: V\rightarrow X'$ be a common resolution such that $q=\phi_W\circ p$. By Lemma \ref{lem: hl23 2.6},
$$p^*K_{\Aa_W}=q^*K_{\Aa'}+F$$
for some $F\geq 0$ that is exceptional over $X'$. Then we have
$$p^*h^*K_{\Aa}=q^*K_{\Aa'}+F-p^*E,$$
so
$$p^*E-F\sim_{\Rr,X}q^*K_{\Aa'}$$
is nef$/X$. Since $h_*p_*(F-p^*E)=h_*p_*F\geq 0$, by the negativity lemma, $F\geq p^*E$. Thus $a(D,\Aa)\leq a(D,\Aa')$ for any prime divisor $D$ over $X$.
\end{proof}

\noindent\textit{Proof of Lemma \ref{lem: g-pair version bir12 2.8} continued}. Let $t$ be the parameter of $\Aa$. First we prove the bs-weak lc model case. Let $\Aa'/U$ be a bs-weak lc model of $\Aa_W/U$ with ambient variety $X'$ and induced birational map $\phi_W: W\dashrightarrow X'$.  By Claim \ref{claim: log smooth model log discrepancy compare}, we only need to show that $\Aa'/U$ is a log birational model of $\Aa/U$. Let $\phi: X\dashrightarrow X'$ be the induced morphism and let
$$\tilde\Aa':=\left(\phi_*\Aa,\Exc(\phi^{-1})^{\ninv}+(1-t)\Exc(\phi^{-1})^{\inv}\right).$$
Then we only need to show that $\Aa'=\tilde\Aa'$. Since  $\Aa'/U$ is a bs-weak lc model of $\Aa_W/U$, we have
$$\Aa'=\left((\phi_W)_*\Aa_W,\Exc(\phi_W^{-1})^{\ninv}+(1-t)\Exc(\phi_W^{-1})^{\inv}\right).$$ 
Let $D$ be a prime divisor on $X'$. Let $\tilde B'$ and $B'$ be the boundaries of $\tilde\Aa'$ and $\Aa'$ respectively. There are three cases:

\medskip

\noindent\textbf{Case 1}. $D$ is not exceptional over $X$. In this case,
    $$-\mult_D\tilde B'=a(D,\tilde\Aa')=a(D,\Aa)=a(D,\Aa_W)=a(D,\Aa')=-\mult_DB',$$
so $\mult_DB'=\mult_D\tilde B'$.

\medskip

\noindent\textbf{Case 2}. $D$ is exceptional over $W$. In this case, $D$ is a component of $\Exc(\phi_W^{-1})$ and a component of $\Exc(\phi^{-1})$, hence
$$\mult_D\tilde B'=t\epsilon_{\Ff}(D)+(1-t)=\mult_DB'.$$
\medskip

\noindent\textbf{Case 3}. $D$ is exceptional over $X$ but not exceptional over $W$. In this case,
$$-\mult_DB'=a(D,\Aa')=a(D,\Aa_W).$$
Since $E\geq 0$,
$a(D,\Aa_W)\leq a(D,\Aa).$
By Claim \ref{claim: log smooth model log discrepancy compare}, 
$a(D,\Aa)\leq a(D,\Aa').$
Thus
$$-\mult_DB'=a(D,\Aa)=a(D,\Aa')=a(D,\Aa_W).$$
By Definition \ref{defn: log smooth models}(4), $a(D,\Aa)=-t\epsilon_{\Ff}(D)-(1-t),$
which implies that
$$\mult_DB'=t\epsilon_{\Ff}(D)+(1-t)=\mult_D\left(\Exc(\phi^{-1})^{\ninv}+(1-t)\Exc(\phi^{-1})^{\inv}\right)=\mult_D\tilde B'.$$

Thus $B'=\tilde B'$, so $\Aa'/U$ is a log birational model of $\Aa/U$, and we are done for the bs-weak lc model case.

Next we prove the bs-minimal model case. Suppose that $\Aa'/U$ be a bs-minimal model of $\Aa_W/U$. For any prime divisor $D$ on $X$ which is exceptional over $X'$, $h^{-1}_*D$ is a prime divisor on $W$ which is exceptional over $X'$. Thus
$$a(D,\Aa)=a(D,\Aa_W)<a(D,\Aa').$$
The bs-minimal model case immediately follows from the bs-weak lc model case. 

The bs-semi-ample model, bs-good minimal model, log minimal model, and good log minimal model cases follow immediately from the bs-weak lc model and the bs-minimal model cases.
\end{proof}

\begin{lem}\label{lem: foliation lsm has lmm}
Let $\Aa/U$ be an lc algebraically integrable adjoint foliated structure and $\Aa'/U$ a bs-weak lc model of $\Aa/U$. Let $\Aa_W$ be a foliated log smooth model of $\Aa$. Let $X,X',W$ be the ambient varieties of $\Aa,\Aa',\Aa_W$ respectively, and assume that the induced birational map $\phi_W: W\dashrightarrow X'$ is a morphism. 

Then we may run a $K_{\Aa_W}$-MMP$/X'$ with scaling of an ample$/X'$ $\Rr$-divisor which terminates with a good minimal model $\Aa_Y/X'$ of $\Aa_W/X'$ such that $$K_{\Aa_Y}=q^*K_{\Aa'}.$$
where $Y$ is the ambient variety of $\Aa_Y$ and $q: Y\rightarrow X'$ is the induced morphism. In particular, $\Aa_Y/U$ is a log minimal model of $\Aa_W/U$.
\end{lem}
\begin{proof}
Let $t$ be the parameter of $\Aa$. If $t=1$ then we are done by \cite[Lemma 4.13]{LMX24}, so we may assume that $t<1$.

Let $h: W\rightarrow X$ be the induced birational morphism. We have
    $$K_{\Aa_W}=h^*K_{\Aa}+E$$
for some $E\geq 0$ that is exceptional$/X$. By Lemma \ref{lem: hl23 2.6}, we have
$$h^*K_{\Aa}=\phi_W^*K_{\Aa'}+F$$
where $F\geq 0$ is exceptional$/X'$. Thus
$$K_{\Aa_W}\sim_{\mathbb R,X'}F+E.$$
\begin{claim}\label{claim: wglc to lmm E exceptional}
$E$ is exceptional$/X'$.
\end{claim}
\begin{proof}
Let $D$ be a component of $E$. By Definition \ref{defn: log smooth models}(2), $a(D,\Aa)>-t\epsilon_{\Ff}(D)-(1-t)$ and $D$ is exceptional$/X$. 

Assume that $D$ is not exceptional over $X'$. Since $\Aa'/U$ is a log birational model of $\Aa/U$ and $\Aa$ is lc, $a(D,\Aa')=-t\epsilon_{\Ff}(D)-(1-t)$. Since $F\geq 0$, $a(D,\Aa)\leq a(D,\Aa')$. Thus $a(D,\Aa)=-t\epsilon_{\Ff}(E)-(1-t)$, hence $D$ is not a component of $E$, a contradiction.
\end{proof}
\noindent\textit{Proof of Lemma \ref{lem: foliation lsm has lmm} continued}. By Claim \ref{claim: wglc to lmm E exceptional}, $F+E$ is exceptional over $X'$. By Theorem \ref{thm: very exceptional mmp}, we may run a $K_{\Aa_W}$-MMP$/X'$ with scaling of an ample$/X'$ divisor, which terminates with a good minimal model $\Aa_Y/X'$ of $\Aa_W/X'$ such that $K_{\Aa_Y}\sim_{\Rr,X'}0$. Since $t<1$, by \cite[Proposition 3.5]{Cas+24}, $\Aa_Y$ is $\Qq$-factorial qdlt, and $a(D,\Aa_W)<a(D,\Aa_Y)$ for any prime divisor $D$ on $W$ that is exceptional$/Y$. By the negativity lemma, 
$$K_{\Aa_Y}=q^*K_{\Aa'}.$$ 
The lemma follows.
\end{proof}

\begin{lem}\label{lem: g-pair weak glc imply lmm}
Let $\Aa/U$ be an lc algebraically integrable adjoint foliated structure. If $\Aa/U$ has a bs-weak lc model (resp. bs-semi-ample model), then $\Aa/U$ has a log minimal model (resp. good log minimal model).
\end{lem}
\begin{proof}
By Lemma \ref{lem: g-pair version bir12 2.7} we only need to prove the bs-weak lc model case. The lemma follows immediately from Lemmas \ref{lem: g-pair version bir12 2.8} and \ref{lem: foliation lsm has lmm}.
\end{proof}

\begin{lem}\label{lem: same weak glc model under pullback}
Let $\Aa/U$ be an lc algebraically integrable adjoint foliated structure with ambient variety $X$. Let $f: Y\rightarrow X$ be a birational morphism, $E\geq 0$ an $f$-exceptional $\mathbb R$-divisor, and $\Aa_Y:=(f^*\Aa,E)$. Assume that $\Aa_Y$ is lc. Then:
\begin{enumerate}
    \item Any bs-weak lc model of $\Aa/U$ is a bs-weak lc model of $\Aa_Y/U$.
    \item If $\Aa/U$ has a bs-weak lc model (resp. bs-semi-ample model), then $\Aa_Y/U$ has a log minimal model (resp. good log minimal model).
\end{enumerate}
\end{lem}
\begin{proof}
(1) Let $\Aa'/U$ be a bs-weak lc model of $\Aa/U$, $\phi: X\dashrightarrow X'$ the induced birational map, and $\phi_Y:=\phi\circ f$. Let $p: W\rightarrow Y$ and $q: W\rightarrow X'$ be a resolution of indeterminacy, and let $h:=f\circ p$. By Lemma \ref{lem: hl23 2.6},
$$h^*K_{\Aa}=q^*K_{\Aa'}+F$$
for some $F\geq 0$ that is exceptional over $X'$. Thus 
$$p^*K_{\Aa_Y}=q^*K_{\Aa'}+p^*E+F.$$
Thus $a(D,\Aa_Y)\leq a(D,\Aa')$ for any prime divisor $D$ over $X'$. In particular, if $a(D,\Aa')=-\epsilon_{\Ff}(D)$, then $a(D,\Aa_Y)=-\epsilon_{\Ff}(D)$.

Let $t$ be the parameter of $\Aa$. Since $\Aa'/U$ is a log birational model of $\Aa/U$ and $\Aa$ is lc, we have
$$\Aa'=\left(\phi_*\Aa,\Exc(\phi^{-1})^{\ninv}+(1-t)\Exc(\phi^{-1})^{\inv}\right).$$
Let 
$$\tilde\Aa'=\left((\phi_Y)_*\Aa_Y,\Exc(\phi_Y^{-1})^{\ninv}+(1-t)\Exc(\phi_Y^{-1})^{\inv}\right).$$ 
Let $B'$ and $\tilde B'$ be the boundaries of $\Aa'$ and $\tilde\Aa'$ respectively. For any prime divisor $D$ on $X'$, there are two cases:

\medskip

\noindent\textbf{Case 1}. $D$ is not exceptional over $X$. In this case,
$$-\mult_DB'=a(D,\Aa')=a(D,\Aa)=a(D,\Aa_Y)=a(D,\Aa')=-\mult_D\tilde B',$$
so $\mult_DB'=\mult_D\tilde B'$.

\medskip

\noindent\textbf{Case 2}. $D$ is exceptional over $X$. In this case, 
$$a(D,\Aa')=-\mult_DB'=-t\epsilon_{\Ff}(D)-(1-t).$$
Since $a(D,\Aa_Y)\leq a(D,\Aa')$ and $\Aa_Y$ is lc, $a(D,\Aa_Y)=-t\epsilon_{\Ff}(D)-(1-t)$. Therefore, if $D$ is not exceptional over $Y$, then
$$\mult_D\tilde B'=-a(D,\Aa_Y)=\epsilon_{\Ff}(D)=\mult_DB',$$
and if $D$ is exceptional over $Y$, then
$$\mult_D\tilde B'=\mult_D\left(\Exc(\phi_Y^{-1})^{\ninv}+(1-t)\Exc(\phi_Y^{-1})^{\inv}\right)=t\epsilon_{\Ff}(D)+(1-t)=\mult_DB'.$$
Thus $B'=\tilde B'$, hence $\Aa'/U$ is a log birational model of $\Aa_Y/U$. Since $K_{\Aa'}$ is nef$/U$, and
$a(D,\Aa_Y)\leq a(D,\Aa')$ for any prime divisor $D$ over $X'$, $\Aa'/U$ is a bs-weak lc model of $\Aa_Y/U$, and we get (1).

(2) It follows from (1), Lemma \ref{lem: g-pair weak glc imply lmm}, and Lemma \ref{lem: g-pair version bir12 2.7}.
\end{proof}

\begin{lem}\label{lem: mm preserved under dlt model}
Let $\Aa/U$ be an lc adjoint foliated structure with ambient variety $X$. Let $f: Y\rightarrow X$ be a birational morphism which only extracts lc places of $\Aa$ and $\Aa_Y:=f^*\Aa$. 

Then any bs-weak lc model (resp. bs-minimal model) of $\Aa_Y/U$ is a bs-weak lc model (resp. bs-minimal model) of $\Aa/U$.
\end{lem}
\begin{proof}
    Let $\Aa'/U$ be a bs-weak lc model (resp. bs-minimal model) of $\Aa_Y/U$ and let $X'$ be its ambient variety. Then for any prime divisor $D$ on $X$ that is exceptional$/X'$, we have
    $$a(D,\Aa)=a(D,\Aa_Y)\leq\text{(resp.}<\text{)} a(D,\Aa').$$
    For any prime divisor $D$ on $X'$ that is exceptional$/X$, if $D$ is exceptional$/Y$, then $D$ is an lc place of $\Aa'$, hence an lc place of $\Aa$. If $D$ is not exceptional$/Y$, then $\Center_YD$ is a prime divisor that is contracted by $f$, so $D$ is an lc place of $\Aa$. Therefore, $\Aa'/U$ is a log birational model of $\Aa/U$. The lemma follows.
\end{proof}

\begin{lem}\label{lem: minimal model same after running mmp}
    Let $\Aa/U$ be an adjoint foliated structure with ambient variety $X$ and let $\phi: X\dashrightarrow X'$ be a $K_{\Aa}$-negative map$/U$ with $\Aa':=\phi_*\Aa$. Then:
    \begin{enumerate}
        \item Any (bs-)minimal model of $\Aa/U$ is a (bs-)minimal model of $\Aa'/U$.
        \item Any (bs-)minimal model of $\Aa'/U$ is a (bs-)minimal model of $\Aa/U$.
    \end{enumerate}
\end{lem}
\begin{proof}
(1) Let $\Aa''/U$ be a bs-minimal model of $\Aa/U$ with ambient variety $X''$. 

For any prime divisor $D$ on $X'$ that is exceptional over $X''$, we have
    $$a(D,\Aa')=a(D,\Aa)<a(D,\Aa'').$$
For any prime divisor $D$ on $X''$ that is exceptional$/X'$, if $D$ is not exceptional$/X$, then by \cite[Lemma 2.25]{LMX24}, $\sigma_D(X/U,K_{\Aa})>0$.  
Let $p: W\rightarrow X$, $q: W\rightarrow X''$ be a resolution of indeterminacy of the induced birational map $\psi: X\dashrightarrow X''$. By Lemma \ref{lem: hl23 2.6},
$$p^*K_{\Aa}=q^*K_{\Aa''}+E$$
for some $E\geq 0$ that is exceptional$/X''$, so
$$0<\sigma_D(X/U,K_{\Aa})=\sigma_{p^{-1}_*D}(W/U,q^*K_{\Aa''}+E)=\sigma_{p^{-1}_*D}(E)=0,$$
which is not possible. Therefore, any prime divisor $D$ on $X''$ that is exceptional$/X'$ is also exceptional$/X$, so $\Aa''/U$ is a log birational model of $\Aa'/U$, and if $\psi$ does not extract any divisor, then the induced birational map $X'\dashrightarrow X''$ does not extract any divisor. (1) follows.

(2) Let $\Aa''/U$ be a bs-minimal model of $\Aa'/U$. Then for any prime divisor $D$ on $X$ that is exceptional$/X''$, if $D$ is not exceptional$/X'$, then
$$a(D,\Aa)=a(D,\Aa')<a(D,\Aa''),$$
and if $D$ is exceptional$/X'$, then
$$a(D,\Aa')=a(D,\Aa)<a(D,\Aa'').$$
Moreover, since $\phi_*\Aa=\Aa'$ and $\Aa''/U$ is a log birational model of $\Aa'/U$, $\Aa''/U$ is a log birational model of $\Aa/U$. Hence $\Aa'/U$ is a bs-minimal model of $\Aa/U$. Moreover, if the induced birational map $X'\dashrightarrow X''$ does not extract any divisor, then  the induced birational map $X\dashrightarrow X''$ does not extract any divisor. The lemma follows.
\end{proof}

\begin{lem}\label{lem: scaling basic properties}
Let $\Aa/U$ be an adjoint foliated structure with ambient variety $X$. Let $\phi_i: X_i\dashrightarrow X_{i+1}$, $X_0:=X$ be a sequence of steps of a $K_{\Aa}$-MMP$/U$ with scaling of some $\mathbb R$-divisor $C$. Let $\Aa_i,C_i$ be the images of $\Aa,C$ on $X_i$ respectively and let
$$\lambda_i:=\inf\{s\geq 0\mid K_{\Aa_i}+sC_i\text{ is nef}/U\}$$
be the scaling numbers of this MMP. Assume this MMP does not terminate. Let $\lambda:=\lim_{i\rightarrow+\infty}\lambda_i$. Then for any $i\gg 0$, we have the following.
\begin{enumerate}
    \item $\phi_i$ is a flip.
    \item If $\lambda=0$, then $K_{\Aa_i}$ is movable$/U$.
    \item Assume that $\lambda=0$ and $\Aa/U$ has a bs-minimal model (resp. minimal model) $\Aa'/U$ with ambient variety $X'$. Let $\psi_i: X_i\dashrightarrow X'$ be the induced birational map. Then $\psi_i$ does not contract any divisor (resp. $\psi_i$ is small).
\end{enumerate}
\end{lem}
\begin{proof}
(1) Let $h: Y\rightarrow X$ be a resolution of $X$, $n$ the number of $h$-exceptional prime divisors, and $\rho:=\dim N^1(Y/U)\otimes{\mathbb R}-n$. 
By \cite[Lemma 1.6]{AHK07}, for any $i$, the induced birational map $X\dashrightarrow X_i$ contracts 
at most $\rho-1$ prime divisors, and (1) follows.

(2) Let $\phi_{i,j}: X_i\dashrightarrow X_j$ be the induced birational maps. By (1), if $i\gg 0$ then
$$K_{\Aa_i}=\lim_{j\rightarrow+\infty}(f_{i,j}^{-1})_*(K_{\Aa_j}+\lambda_jC_j)$$
is movable$/U$. 

(3) Let $p: W\rightarrow X_i$ and $q: W\rightarrow X'$ be a resolution of indeterminacy of $\psi_i$. By Lemma \ref{lem: minimal model same after running mmp}, $\Aa'/U$ is a bs-minimal model (resp. minimal model) of $\Aa_i/U$. By Lemma \ref{lem: hl23 2.6},
$$F:=p^*K_{\Aa_i}-q^*K_{\Aa'}\geq 0$$
is exceptional$/X'$. By (2) and \cite[Lemma 3.4]{LX25},
$$0=N_{\sigma}(K_{\Aa_i}/U)=p_*N_{\sigma}(p^*K_{\Aa_i}/U)=p_*N_{\sigma}((q^*K_{\Aa'}+F)/U)=p_*F.$$
Since $a(D,\Aa_i)<a(D,\Aa')$ for any prime divisor $D$ on $X_i$ that is exceptional$/X'$, $\Supp(p_*F)$ contains all prime divisors on $X_i$ that are exceptional$/X'$, hence $\psi_i$ does not contract any divisor. In particular, if $\Aa'/U$ is a minimal model of $\Aa/U$, then $\psi_i$ is small.
\end{proof}

\section{Existence of minimal models and MMP with scaling}\label{sec: mmp scaling}

The goal of this section is to prove that, for lc algebraically integrable adjoint foliated structures, existence of a minimal model, even in the sense of Birkar-Shokurov, implies the termination of MMP with scaling. This result can be seen as an analogue of \cite[Theorem 1.9]{Bir12}, \cite[Theorem 4.1]{HL22}, \cite[Theorem 4.1]{LT22}, \cite[Theorem 2.19]{TX24} for pairs or generalized pairs. Our proof generally follows the same ideas in these references with minor differences as we do not need to introduce an auxiliary ample divisor whose components span the Neron-Severi cone. The main result of this section is the following:

\begin{thm}\label{thm: eomm implies tof with scaling}
    Let $\Aa(0)/U:=(X,\Ff,B_0^{\ninv}+(1-t_0)B_0^{\inv},\Mm_0,t_0)/U$ and $\Aa(1)/U:=(X,\Ff,B_1^{\ninv}+(1-t_1)B_1^{\inv},\Mm_1,t_1)/U$ be two lc algebraically integrable adjoint foliated structures. Let $\Aa(s):=s\Aa(1)+(1-s)\Aa(0)$ for any $s\in\mathbb R$ and let $C:=K_{\Aa_1}-K_{\Aa_0}$. Assume the following.
    \begin{enumerate}
        \item $B_1\geq B_0$, $\Mm_1-\Mm_0$ is nef$/U$, and $t_1\geq t_0$.
        \item $\phi_i: X_i\dashrightarrow X_{i+1}$ is a sequence of steps of a $K_{\Aa(0)}$-MMP$/U$ with scaling of $C$ with $X_0:=X$, $\Aa_i(s),C_i$ the images of $\Aa(s),C$ on $X_i$ respectively for any $s\in [0,1]$ and any $i$, and
    $$\lambda_i:=\inf\{s\geq 0\mid K_{\Aa_i(s)}\text{ is nef}/U\}$$
    the scaling numbers of this MMP. In particular, we assume that $\lambda_0\leq 1$,  $K_{\Aa(\lambda_0)}$ is nef, and there exists a $K_{\Aa_i(0)}$-negative extremal ray$/U$ that is $K_{\Aa_i(\lambda_i)}$-trivial unless $\lambda_i=0$.
    \item $K_{\Aa_i(\lambda_i)}$ is NQC$/U$ for any $i\gg 0$.
    \item $\Aa(0)/U$ has a KNQC bs-weak lc model $\Aa_Y/U$.
    \end{enumerate}
Then either this MMP terminates or $\lim_{i\rightarrow+\infty}\lambda_i>0$.
\end{thm}

We need to prove several lemmas before proving Theorem \ref{thm: eomm implies tof with scaling}.

\begin{lem}\label{lem: qdlt decrease coefficient}
    Let $\Aa_i/U:=(X,\Ff,B_i^{\ninv}+(1-t_i)B_i^{\inv},\Mm_i,t_i)/U$, $i\in\{1,2\}$ be two algebraically integrable adjoint foliated structures such that $B_2\geq B_1,\Mm_2-\Mm_1$ is nef$/U$, and $t_2\geq t_1$. If $\Aa_2$ is $\mathbb Q$-factorial qdlt, then $\Aa_1$ is  $\mathbb Q$-factorial qdlt.
\end{lem}
\begin{proof}
We first reduce to the case when $t_1=t_2$. If $t_1<t_2$, then $$\Aa_1':=(X,\Ff,B_2^{\ninv}+(1-t_1)B_2^{\inv},\Mm_2,t_1)=\frac{t_1}{t_2}\Aa_2+\frac{t_2-t_1}{t_2}(X,B_2,\Mm_2)$$
is lc, and and any lc place of $\Aa_1'$ is an lc place of $(X,B_2,\Mm_2)$. Moroever, since $\Aa_2$ is qdlt, $t_2<1$ and $(X,B_2,\Mm_2)$ is qdlt. Thus $\Aa_1'$ is qdlt. Possibly replacing $\Aa_2$ with $\Aa'_1$, we may assume that $t_1=t_2$. 

Let $\Nn:=\Mm_2-\Mm_1$. Let $E$ be an lc place of $\Aa_1$. We let $h: Y\rightarrow X$ be a birational morphism such that $\Mm_1,\Mm_2$ descend to $Y$ and $E$ is on $Y$. Since $B_2\geq B_1$, $\Mm_2-\Mm_1$ is nef, and $\Aa_2$ is lc, we have that $E$ is an lc place of $\Aa_2$, $\mult_E(B_2-B_1)=0$, and $\mult_E(h^*\Nn_{X}-\Nn_{Y})=0$. Since $\Aa_2$ is qdlt, $E$ is an lc place of $(X,B_2,\Mm_2)$. We have
$$a(E,X,B_1,\Mm_1)=a(E,X,B_2,\Mm_2)+\mult_E(B_2-B_1)+\mult_E(h^*\Nn_{X}-\Nn_Y)=-1.$$
Thus $E$ is an lc place of $(X,B_1,\Mm_1)$. By definition, $\Aa_1$ is qdlt.
\end{proof}

\begin{lem}\label{lem: trivial mmp}
    Let $\Aa/U=(X,\Ff,B,\Mm,t)/U$ be an lc algebraically integrable adjoint foliated structure such that $X$ is potentially klt and $K_{\Aa}$ is NQC$/U$. Write $K_{\Aa}=\sum a_iM_i$ where each $M_i$ is nef$/U$ Cartier and each $a_i>0$. Let $l:=\frac{2\dim X+1}{\min\{a_i\}}$.
    
    Let $\tilde\Aa/U=(X,\Ff,\tilde B,\tilde\Mm,\tilde t)/U$ be an lc algebraically integrable adjoint foliated structure and let $\Aa(s):=s\bar\Aa+(1-s)\Aa$ for any $s\in [0,1]$. Then for any $s\in\left[0,\frac{1}{l+1}\right]$, any sequence of steps of a $K_{\Aa(s)}$-MMP$/U$ is $K_{\Aa}$-trivial.
\end{lem}
\begin{proof}
We may assume that $s>0$. We have
$$K_{\Aa(s)}=s\left(K_{\tilde\Aa}+\left(\frac{1}{s}-1\right)K_{\Aa}\right)$$
so any sequence of steps of a $K_{\Aa(s)}$-MMP$/U$ is a sequence of steps of a $(K_{\tilde\Aa}+rK_{\Aa})$-MMP$/U$ for some $r\geq l$. Assume that a sequence of steps $\phi: X\dashrightarrow X'$ of a $K_{\Aa(s)}$-MMP$/U$ is $K_{\Aa}$-trivial. Let $\Aa'(a):=\phi_*\Aa(a)$ for any $a\in [0,1]$ and $M_i':=\phi_*M_i$ for each $i$. Then by the contraction theorem, $M_i'$ is nef$/U$ and Cartier for any $i$. Moreover, $\phi$ is also $K_{\Aa(1)}$-negative, so $\Aa'(1)$ is lc.

For any step $X'\xrightarrow{f} T\leftarrow X''$ of a $K_{\Aa'(s)}$-MMP$/U$ where $f$ is the contraction of an extremal ray, let $R$ be the $K_{\Aa'(s)}$-negative extremal ray contraction by $f$. Then $R$ is also $K_{\Aa'(1)}$-negative. By \cite[Theorem 1.3]{Cas+24}, $R$ is spanned by a rational extremal ray $C$ such that $0>K_{\Aa'(1)}\cdot C\geq -2\dim X$. If $K_{\Aa'(0)}\cdot C>0$, then
\begin{align*}
    0>&K_{\Aa'(s)}\cdot C=s\left(K_{\Aa'(1)}+\left(\frac{1}{s}-1\right)K_{\Aa'(0)}\right)\cdot C\geq s\left(-2\dim X+l\sum a_iM_i'\cdot C\right)\\
    &\geq s\left(-2\dim X+\frac{2\dim X+1}{\min\{a_i\}}\cdot\min\{a_i\}\right)>0,
\end{align*}
which is not possible. Thus $f$ is $K_{\Aa'(0)}$-trivial, and we are done by induction on the number of steps of the MMP.
\end{proof}

\begin{lem}[Liftying the MMP]\label{lem: lift mmp afs}
Let $\Aa/U$ be an lc algebraically integrable adjoint foliated structure with ambient variety $X$ and parameter $t$. Let $C$ be an $\mathbb R$-divisor on $X$ and let $\phi_i: X_i\xrightarrow{f_i}T_i\xleftarrow{f_i^+} X_{i+1}$, $X_0:=X$ be a sequence of steps of a $K_{\Aa}$-MMP$/U$ and $\Aa_i,C_i$ the images of $\Aa,C$ on $X_i$ respectively, where $f_i$ is the contraction of an extremal ray. Let $h: X'\rightarrow X$ be a $\mathbb Q$-factorial ACSS modification of $\Aa$ if $t=1$ and a $\mathbb Q$-factorial qdlt modification of $\Aa$ if $t<1$. Let $\Aa':=h^*\Aa$.

Assume that either $X$ is potentially klt, or Theorem \ref{thm: eomm implies tof with scaling} holds when $X$ is $\mathbb Q$-factorial klt.  Then there exists a sequence of birational maps$/U$ $\phi_i': X_i'\dashrightarrow X_{i+1}'$, $X_0':=X'$ and birational morphisms $h_i: X_i'\rightarrow X_i$, $h_0:=h$ satisfying the following. Let $\Aa_i'$ be the image of $\Aa'$ on $X_i'$, then for each $i$:
\begin{enumerate}
    \item $h_i$ is a $\mathbb Q$-factorial ACSS modification of $\Aa_i$ if $t=1$ and is a $\mathbb Q$-factorial qdlt modification of $\Aa_i$ when $t<1$, with $h_i^*\Aa_i=\Aa_i'$.
    \item $h_{i+1}\circ\phi_i'=\phi_i\circ h_i$.
    \item $\phi_i'$ is a $K_{\Aa_i'}$-MMP$/T_i$ that is not the identity morphism, and $\Aa_{i+1}'/T_i$ is a good log minimal model of $\Aa_i'/T_i$.
    \item $\{\phi_i'\}$ is a sequence of steps of a $K_{\Aa'}$-MMP$/U$.
    \item If $\{\phi_i\}$ is a sequence of steps of a $K_{\Aa}$-MMP$/U$ with scaling of $C$ with scaling numbers $$\lambda_i:=\inf\{s\geq 0\mid K_{\Aa_i}+sC_i\text{ is nef}/U\},$$ then:
    \begin{enumerate}
        \item $\phi_i'$ is a sequence of steps of a $K_{\Aa'}$-MMP$/U$ with scaling of $C':=h^*C$ with scaling numbers all equal to $\lambda_i$.
        \item $\{\phi_i'\}$ is a sequence of steps of a $K_{\Aa'}$-MMP$/U$ with scaling of $C'$.
    \end{enumerate} 
\end{enumerate}
\end{lem}
\begin{proof}
If $t=1$, then we are done by \cite[Propositions 8.2, A.45]{LMX24}. Note that \cite{LMX24} further requires that $h$ is ``strict and super" but the same lines of the arguments work in our setting. In the rest of the proof, we assume that $t<1$.

Since (4) follows from (3) 
    and (5.b) follows from (5.a), we only need to prove (1)(2)(3) and (5.a).

Let $n$ be a non-negative integer. We prove the proposition by induction on $n$ under the induction hypothesis that we have already constructed $\Aa_i'/U$ and $h_i$ for any $i\leq n$ and $\phi_{i,Y}$ for any $i\leq n-1$ which satisfy (1)(2)(3)(5). When $n=0$, this follows from our assumption, so we may assume that $n>0$. We need to construct $\phi_{n}',h_{n+1}$, and $\Aa_{n+1}'/U$.

When $X$ is potentially klt, we let $(X,\Delta)$ be a klt pair, $$\Aa(\epsilon):=\epsilon(X,\Delta)+(1-\epsilon)\Aa$$
for any $\epsilon\in [0,1]$, and $\Aa_i(\epsilon)$ the image of $\Aa(\epsilon)$ on $X_i$ for any $i$. Then $\Aa(\epsilon)$ is klt for any $\epsilon\in (0,1]$, and the induced birational map $X\dashrightarrow X_n$ is also a sequence of steps of a $K_{\Aa(\epsilon)}$-MMP for any $0\leq\epsilon\ll 1$. Thus $\Aa_n(\epsilon)$ is klt for any $0<\epsilon\ll 1$. Moreover, $\phi_{n}$ is a step of a $K_{\Aa_n(\epsilon)}$-MMP$/U$ for any $0<\epsilon\ll 1$. Let $\Aa_n'(\epsilon):=h_n^*\Aa_n(\epsilon)$ for any $\epsilon\in [0,1]$. Then $\Aa_n'(\epsilon)$ is $\mathbb Q$-factorial sub-klt for any $0<\epsilon\ll 1$. Since $h_n$ only extracts lc places of $\Aa_n$ and $t<1$, $\Aa_n'(\epsilon)$ is $\mathbb Q$-factorial klt for any $0<\epsilon\ll 1$.

We let $H_n$ be a supporting function of the extremal ray$/U$ contracted by $f_n$ and let  
$$L_n:=H_n-K_{\Aa_n},$$
such that $L_n=\lambda_nC_n$ if we are running an MMP$/U$ with scaling of $C_n$. Then $L_n$ is ample$/T_n$. 

We run a $K_{\Aa_n'}$-MMP$/T_n$ with scaling of an ample divisor $A_n$. Note that $\Aa_{n+1}/T_n$ is a good minimal model of $\Aa_n'/T_n$. When $X$ is potentially klt, this MMP is also a $K_{\Aa_n'(\epsilon)}$-MMP$/T_n$ with scaling of $\mu A_n$ for some $0<\epsilon\ll 1$. Thus by \cite[Theorem 2.1.1]{Cas+25}, and by Theorem \ref{thm: eomm implies tof with scaling} when $X$ is not potentially klt, this MMP terminates with a good minimal model $\Aa_{n+1}'(\epsilon)/T_n$ of $\Aa_n'(\epsilon)/T_n$ with ambient variety $X_{n+1}$ and induced birational map $\phi_n': X_n'\dashrightarrow X_{n+1}'$. 

Let $\Aa_{n+1}':=(\phi_n')_*\Aa_n'$. By Lemma \ref{lem: numerical equivalence model}, $\Aa_{n+1}'/T_n$ is a good log minimal model of $\Aa_n'/T_n$. This implies (3).

Since $\phi_n$ is the ample model$/T_n$ of $K_{\Aa_n}$, there exists an induced birational morphism $h_{n+1}: X_{n+1}'\rightarrow X_n$ such that $h_{n+1}\circ\phi_n'=\phi_n\circ h_n$. This implies (2). Moreover, $K_{\Aa_{n+1}'}=h_{n+1}^*K_{\Aa_{n+1}}$, so $h_{n+1}$ is a $\mathbb Q$-factorial qdlt modification of $\Aa_i$. This implies (1).

Moreover, this MMP is also an MMP$/U$ with scaling of $L_n$ as $K_{\Aa_n'}+h_n^*L_n\sim_{\mathbb R,T_n}0$ and $K_{\Aa_n'}+h_n^*L_n=h^*H_n$ is nef$/U$. This implies (5.a) and we are done.
\end{proof}

\begin{rem}\label{rem: may assume 4.1}
We are going to prove Theorem \ref{thm: eomm implies tof with scaling} in the following way: we are going to use Lemma \ref{lem: lift mmp afs} for the $\mathbb Q$-factorial klt case to prove Theorem \ref{thm: eomm implies tof with scaling} when $X$ is $\mathbb Q$-factorial klt. This in turn imply Lemma \ref{lem: lift mmp afs} unconditionally. Using the same lines of the proof, we can prove the general case of Theorem \ref{thm: eomm implies tof with scaling}. 
\end{rem}

\begin{proof}[Proof of Theorem \ref{thm: eomm implies tof with scaling}]
    We prove the Theorem contradiction. By Lemma \ref{lem: g-pair weak glc imply lmm} we may assume that $\Aa_Y/U$ is a $\mathbb Q$-factorial log minimal model of $\Aa(0)/U$.
    We may assume that this MMP does not terminates and $\lim_{i\rightarrow+\infty}\lambda_i=0$. Let $\varphi: X\dashrightarrow Y$ be the induced birational map.

    By Remark \ref{rem: may assume 4.1}, we may assume that one of the following conditions holds:
    \begin{itemize}
        \item $X$ is $\mathbb Q$-factorial klt.
        \item The $\mathbb Q$-factorial klt case of Theorem \ref{thm: eomm implies tof with scaling} holds.
    \end{itemize}

    \medskip

\noindent\textbf{Step 1.} Write $\Aa_i(s):=(X_i,\Ff_i,B_i(s),\Mm(s),t(s))$ for any $i$ and any $s\in [0,1]$. By Lemma \ref{lem: scaling basic properties}, possibly truncating the MMP, we may assume that 
    \begin{itemize}
    \item[(i)] $\phi_i$ is a flip for every $i$ and $K_{\Aa(0)}$ is movable$/U$.
    \end{itemize}
 Possibly replacing $\Aa(s)$ with $\Aa_{i+1}(s)$ for any $s\in [0,1]$ and some $i+1$ such that $\lambda_i>\lambda_{i+1}$, and then replacing $\Aa(1)$ with $\Aa(\lambda_1)$,
 we may assume that
    \begin{itemize}
        \item[(ii)] $\lambda_1=1$, and
        \item[(iii)] $\Aa(1+\epsilon)$ is lc for some $\epsilon>0$. In particular, if $t_0<1$, then $t_1<1$.
\end{itemize}
Write $\Aa_Y(0):=(Y,\Ff_Y,B_Y,\Mm,t_0)/U$.
Let $p: W\rightarrow X$ and $q: W\rightarrow Y$ be a resolution of indeterminacy of $\varphi$ so that $p$ is a foliated log resolution of $(X,\Ff,\Supp B_0\cup\Supp B_1,\Mm)$ and $q$ is a foliated log resolution of $(Y,\Ff,\Supp B_Y,\Mm)$. By Lemma \ref{lem: scaling basic properties}, $\varphi$ does not contract any divisor.

\medskip

\noindent\textbf{Step 2.} In this step we construct a model $Y'$. Let 
$$\Aa_{W}(s):=\left(p^{-1}_*\Aa(s),\Exc(p)^{\ninv}+(1-st_1-(1-s)t_0)\Exc(p)^{\inv}\right)$$
Then $\Aa_W(s)$ is a foliated log smooth model of $\Aa(s)$ for any $s\in [0,1+\epsilon]$. By Lemma \ref{lem: foliation lsm has lmm}, we may run a $K_{\Aa_W(0)}$-MMP$/Y$ with scaling of an ample divisor which terminates with a good log minimal model $\Aa_{Y'}(0)/Y$ of $\Aa_W(0)/Y$ such that $K_{\Aa_{Y'}(0)}\sim_{\mathbb R,Y}0$. Let $Y'$ be the ambient variety of $\Aa_{Y'}(0)$, $\psi: W\dashrightarrow Y'$ the induced birational map, and $\Aa_{Y'}(s):=\psi_*\Aa_W(s)$ for any $s$.

\medskip

\noindent\textbf{Step 3.} In this step we construct a model $X'$. Let $F(s):=K_{\Aa_W(s)}-p^*K_{\Aa(s)}$ for any $s$. Since $K_{\Aa_W(1)}\sim_{\mathbb R,X}F(1)\geq 0$
is exceptional$/X$, by Theorem \ref{thm: very exceptional mmp}, we may run a $K_{\Aa_W(1)}$-MMP$/X$ which terminates with a good log minimal model $\Aa'(1)/X$ of $\Aa(1)/X$ with ambient variety $X'$ and associated with birational map $\xi: W\dashrightarrow X'$ such that $K_{\Aa'(1)}\sim_{\mathbb R,X}0$. Let $\Aa'(s):=\xi_*\Aa_W(s)$ for any $s$. We have $K_{\Aa'(1)}=h^*K_{\Aa(1)}$ where $h: X'\rightarrow X$ is the induced birational morphism.

By (iii), $\Aa'(1)$ is $\mathbb Q$-factorial qdlt when $t_0<1$ and is $\mathbb Q$-factorial ACSS when $t_0=1$. By Lemma \ref{lem: qdlt decrease coefficient}, for any $s\in [0,1]$, $\Aa'(s)$ is $\mathbb Q$-factorial qdlt when $t_0<1$ and is $\mathbb Q$-factorial ACSS when $t_0=1$. 

The divisors contracted by $\xi$ are exactly the prime divisors $W$ that are irreducible components of $F(1)$, which are exactly prime divisors on $W$ that are not lc places of $\Aa(1)$. Since $\Aa(1+\epsilon)$ is lc and
$$\Aa(1)=\frac{1}{1+\epsilon}\Aa(1+\epsilon)+\frac{\epsilon}{1+\epsilon}\Aa(0),$$
any lc place of $\Aa(1)$ is an lc place of $\Aa(0)$. By (1) and \cite[Proposition 3.3]{Cas+24}, any lc place of $\Aa(0)$ is an lc place of $\Aa(1)$.

Thus the divisors contracted by $\xi$ are exactly the prime divisors on $W$ that are not lc places of $\Aa(0)$. In particular, $h$ is a $\mathbb Q$-factorial ACSS modification of $\Aa(0)$ if $t_0=1$ and is a $\mathbb Q$-factorial qdlt modification of $\Aa(0)$ if $t_0<1$.

\medskip

\noindent\textbf{Step 4.} We show that $\Aa_{Y'}(0)/U$ is a minimal model of $\Aa'(0)/U$. 

First we show that the induced birational map $\varphi': X'\dashrightarrow Y'$ does not extract any divisor. Suppose that there exists a prime divisor $D$ on $Y'$ that is exceptional$/X'$. Then $D$ is exceptional$/X$, so $D$ is an lc place of $\Aa(0)$. Since $\Center_WD$ is a divisor and $\xi$ only contracts prime divisors that are not lc places of $\Aa(0)$, $D$ is not contracted by $\xi$, a contradiction. 

Now for any prime divisor $D$ on $X'$ that is exceptional$/Y'$, we have
$$a(D,\Aa'(0))=a(D,\Aa_W(0))<a(D,\Aa_{Y'}(0)).$$
Thus $\Aa_{Y'}(0)/U$ is a minimal model of $\Aa'(0)/U$.

\medskip

\noindent\textbf{Step 5.} In this step we construct models $X_i'$.

By Lemma \ref{lem: lift mmp afs}, there exists a sequence of steps of a $K_{\Aa'(0)}$-MMP$/U$ with scaling of $C':=h^*C=K_{\Aa'(1)}-K_{\Aa'(0)}$, $\phi_i': X_i'\dashrightarrow X_{i+1}'$, $X_0':=X'$ satisfying the following. Let $\Aa'_i(s),C_i'$ be the images of $\Aa'(s)$, $C'$ on $X_i'$ respectively for any $i,s$ and
$$\lambda_i':=\inf\{s\geq 0\mid K_{\Aa_i'(s)}\text{ is nef}/U\}$$
the scaling numbers. Then there exist a strictly increasing sequence $n_i$ so that the induced birational maps $h_i: X'_{n_i}\rightarrow X_i$ are morphisms, each $h_i$ is a $\mathbb Q$-factorial ACSS modification of $\Aa_i(0)$ if $t=1$, a $\mathbb Q$-factorial qdlt modification of $\Aa_i(0)$ if $t<1$, $\Aa_{n_i}'(0)=h_i^*\Aa_i(0)$, and $\lambda_j'=\lambda_i$ when $n_i\leq j\leq n_{i+1}-1$. In particular, $\lim_{i\rightarrow+\infty}\lambda_i'=0$.

\medskip

\noindent\textbf{Step 6.} In this step we construct models $X(\delta,\delta')$ for $0<\delta,\delta'\ll 1$.  Let $L$ be an ample divisor on $W$ and $\Ll:=\overline{L}$. By our construction, $\psi$ is a sequence of steps of a $(K_{\Aa_W(\delta)}+\delta'L)$-MMP$/Y$ for any $0\leq \delta,\delta'\ll 1$.

In particular, for any $0\leq\delta,\delta'\ll 1$, $(\Aa_{Y'}(\delta),\delta'\Ll)$ is $\mathbb Q$-factorial qdlt if $t<1$, and is $\mathbb Q$-factorial ACSS if $t=1$. Since $K_{\Aa_{Y'}(0)}$ is NQC$/U$, by Lemma \ref{lem: trivial mmp}, for any $0<\delta,\delta'\ll 1$, any sequence of steps of a $(K_{\Aa_{Y'}(\delta)}+\delta'\Ll_{Y'})$-MMP$/U$ is $K_{\Aa_{Y'}(0)}$-trivial. By \cite[Lemma 7.2]{Cas+25} and \cite[Lemma 16.1.1]{CHLX23}, for any $0<\delta,\delta'\ll 1$, there exists an lc algebraically integrable adjoint foliated structure $\Aa_{Y'}(\delta,\delta')$ with ambient variety $Y'$, and an ample$/U$ $\mathbb R$-divisor $H(\delta,\delta')$ on $Y'$, such that
$$K_{\Aa_{Y'}(\delta,\delta')}+H(\delta,\delta')\sim_{\mathbb R,U}K_{\Aa_{Y'}(\delta)}+\delta'\Ll_Y.$$
By \cite[Theorem 8.1]{Cas+25}, we may run a $(K_{\Aa_{Y'}(\delta,\delta')}+H(\delta,\delta'))$-MMP$/U$ which terminates with a good minimal model $\Aa(\delta,\delta')/U$ of $(\Aa_{Y'}(\delta,\delta'),\overline{H(\delta,\delta')})/U$ with ambient variety $X(\delta,\delta')$ and induced birational map $Y'\dashrightarrow X(\delta,\delta')$. Let $\Aa_{X(\delta,\delta')}(s)$ be the image of $\tilde\Aa_{Y'}(s,\delta')$ on $X(\delta,\delta')$ for any $s$. Then
$$K_{\Aa(\delta,\delta')}\sim_{\mathbb R,U}K_{\Aa_{X(\delta,\delta')}(\delta)}+\delta'\Ll_{X(\delta,\delta')}$$
is semi-ample$/U$ and
$$K_{\Aa_{X(\delta,\delta')}(0)}$$
is nef$/U$ as $Y'\dashrightarrow X(\delta,\delta')$ if $K_{\Aa_{Y'}(0)}$-trivial. Moreover, since $K_{\Aa_{Y'}(\delta)}$ and $\delta'L_{Y'}$ are movable$/U$, $K_{\Aa_{Y'}(\delta,\delta')}+H(\delta,\delta')$ is movable$/U$, so the induced birational map $Y'\dashrightarrow X(\delta,\delta')$ is small.

\medskip

\noindent\textbf{Step 7.} By our construction, $\xi$ is a sequence of steps of a $(K_{\Aa_W(1)}+\delta L)$-MMP$/X$ for any $0<\delta\ll 1$. Since $\Aa_{Y'}(0)/U$ is a minimal model of $\Aa'(0)/U$, by Lemma \ref{lem: minimal model same after running mmp}, $\Aa_{Y'}(0)/U$ is a minimal model of $\Aa_i'(0)/U$ for any $i$. By Lemma \ref{lem: scaling basic properties} and our assumption, we may let $i\gg 0$ be a positive integer such that 
\begin{itemize}
    \item $1\gg\lambda_{i}>\lambda_{i+1}$,
    \item $K_{\Aa_{n_i}'(0)}$ is movable$/U$,
    \item the induced birational map $f: X_{n_i}'\dashrightarrow Y'$ is small,
    \item $K_{\Aa_{i}(\lambda_{i})}$ is NQC$/U$, and
    \item $K_{\Aa_X(\lambda_i,\delta')}+\delta'\Ll_{X(\lambda_i,\delta')}$ is semi-ample$/U$, $K_{\Aa_X(\lambda_i,0)}$ is nef$/U$, and the induced birational map $Y'\dashrightarrow X(\lambda_i,\delta')$ is small for any $0<\delta'\ll 1$.
\end{itemize}
By our construction, $h_{i}^*K_{\Aa_{i}(\lambda_{i})}=K_{\Aa'_{n_{i}}(\lambda_{i})}$ is NQC$/U$. Since the induced birational map $X'\dashrightarrow X_{n_i}$ is a $K_{\Aa'(\lambda_i)}$-MMP$/U$, $X'\dashrightarrow X_{n_i}$ is also a sequence of steps of a $(K_{\Aa'(\lambda_i)}+\delta\Ll_{X'})$-MMP$/U$ for any $0\leq\delta\ll 1$.

By \cite[Lemma 7.2]{Cas+25} and \cite[Lemma 16.1.1]{CHLX23}, for any $0<\delta\ll 1$, there exists an adjoint foliated structure $(\Aa'(1),\Nn(\delta))/U$ for some nef$/U$ $\bb$-divisor $\Nn(\delta)$ and an ample $\mathbb R$-divisor $P(\delta)$ on $X'$, such that
$$K_{\Aa'(1)}+\Nn(\delta)_{X'}+P(\delta)\sim_{\mathbb R,U}K_{\Aa'(1)}+\delta\Ll_{X'},$$
$(\Aa'(1),\Nn(\delta))$ is $\mathbb Q$-factorial ACSS if $t=1$ and $\mathbb Q$-factorial qdlt if $t<1$.
Let $\Pp(\delta):=\overline{P(\delta)}$. For any $0<\delta\ll 1$, by \cite[Lemma 7.2]{Cas+25} and \cite[Lemma 16.1.1]{CHLX23} again, there exists an lc adjoint foliated structure $\mathfrak{B}(\delta)/U$ with ambient variety $X_{n_i}'$ and an ample$/U$ $\mathbb R$-divisor $H(\delta)$ on $X_{n_i}'$, such that
$$K_{\mathfrak{B}(\delta)}+H(\delta)\sim_{\mathbb R,U}K_{\Aa_{n_i}'(\lambda_i)}+\Nn(\delta)_{X_{n_i}'}+\Pp(\delta)_{X_{n_i}'}\sim_{\mathbb R,U}K_{\Aa_{n_i}'(\lambda_i)}+\delta\Ll_{X_{n_i}'}.$$
By Lemma \ref{lem: trivial mmp}, for any $0<\delta\ll 1$, any sequence of steps of a $(K_{\Aa_{n_i}'(\lambda_i)}+\delta\Ll_{X_{n_i}'})$-MMP$/U$ is $K_{\Aa_{n_i}'(\lambda_i)}$-trivial. 

Pick $0<\tau\ll 1$ that is general in $\mathbb R/\mathbb Q$. By \cite[Theorem 8.1]{Cas+25}, we may run a $(K_{\mathfrak{B}(\tau)}+H(\tau))$-MMP$/U$ which terminates with a good minimal model $(\mathfrak{B}_T(\tau),\overline{H(\tau)})/U$ of  $(\mathfrak{B}(\tau),\overline{H(\tau)})/U$ associated with birational map $X_{n_i}'\dashrightarrow T$. Let $\Aa_T(s)$ be the image of $\Aa_{n_i}'(s)$ on $T$ for any $s$, then
$$K_{\Aa_T(\lambda_i)}+\tau\Ll_T$$
is semi-ample$/U$ and
$$K_{\Aa_T(\lambda_i)}$$
is nef$/U$. Moreover, since $K_{\Aa_{n_i}'(\lambda_i)}$ and $L$ are movable$/U$, the induced birational map  $X_{n_i}'\dashrightarrow T$ is small.

Moreover, we have that 
$$K_{\Aa_{X(\lambda_i,\tau)}(\lambda_i)}+\tau\Ll_{X(\lambda_i,\tau)}$$
is semi-ample$/U$ and
$$K_{\Aa_{X(\lambda_i,\tau)}(0)}$$ is nef$/U$, and the induced birational map $Y'\dashrightarrow X(\lambda_i,\tau)$ is small. Thus the induced birational map $\rho: T\dashrightarrow X(\delta,\tau)$ is small.

\medskip

\noindent\textbf{Step 8.} Since $\rho$ is small, $K_{\Aa_{X(\lambda_i,\tau)}(\lambda_i)}+\tau\Ll_{X(\lambda_i,\tau)}$ and $K_{\Aa_T(\lambda_i)}+\tau\Ll_T$ are crepant. Since $\tau$ is general in $\mathbb R/\mathbb Q$, $K_{\Aa_{X(\lambda_i,\tau)}(\lambda_i)}$ and $K_{\Aa_T(\lambda_i)}$ are crepant, hence $K_{\Aa_{X(\lambda_i,\tau)}(\lambda_i)}$ is nef$/U$. Since $K_{\Aa_{X(\lambda_i,\tau)}(0)}$ is nef$/U$ and $\lambda_i>\lambda_{i+1}>0$, $K_{\Aa_{X(\lambda_i,\tau)}(\lambda_{i+1})}$ is nef$/U$.

We have that $K_{\Aa_{n_{i+1}}'(\lambda_i)}$ and $K_{\Aa_{n_i}'(\lambda_i)}$ are crepant, $K_{\Aa_{n_i}'(\lambda_i)}$ and $K_{\Aa_T(\lambda_i)}$ are crepant, and $K_{\Aa_T(\lambda_i)}$ and $K_{\Aa_{X(\lambda_i,\tau)}(\lambda_i)}$ are crepant. Thus $K_{\Aa_{n_{i+1}}'(\lambda_i)}$ and $K_{\Aa_{X(\lambda_i,\tau)}(\lambda_i)}$ are crepant. Since $K_{\Aa_{n_{i+1}}'(\lambda_{i+1})}$ and $K_{\Aa_{X(\lambda_i,\tau)}(\lambda_{i+1})}$ are nef$/U$ and the induced birational map $X_{n_{i+1}}'\dashrightarrow X(\lambda_i,\tau)$ is small, by applying the negativity lemma twice, $K_{\Aa_{n_{i+1}}'(\lambda_{i+1})}$ and $K_{\Aa_{X(\lambda_i,\tau)}(\lambda_{i+1})}$ are crepant. Since $\lambda_i>\lambda_{i+1}$, $K_{\Aa_{n_{i+1}}'(s)}$ and $K_{\Aa_{X(\lambda_i,\tau)}(s)}$ are crepant for any $s$. In particular, $K_{\Aa_{n_{i+1}}'(0)}$ and $K_{\Aa_{X(\lambda_i,\tau)}(0)}$ are crepant, hence $K_{\Aa_{n_{i+1}}'(0)}$ is nef. Thus the MMP $\{\phi_i'\}$ terminates, hence the MMP $\{\phi_i\}$ terminates, a contradiction. We conclude the proof.
\end{proof}

\begin{thm}\label{thm: eo bswlcm implies eomm}
    Let $\Aa/U$ be an lc algebraically integrable adjoint foliated structure such that $X$ is potentially klt. Assume that $\Aa/U$ has a KNQC bs-weak lc model $\Aa_Y/U$. Then any sequence of steps of a $K_{\Aa}$-MMP$/U$ with scaling of an ample$/U$ $\mathbb R$-divisor terminates with a minimal model of $\Aa/U$.
\end{thm}
\begin{proof}
Let $\phi_i:\Aa_i\dashrightarrow\Aa_{i+1}$, $\Aa_0:=\Aa$ be a sequence of steps of a  $K_{\Aa}$-MMP$/U$ with scaling of an ample$/U$ $\mathbb R$-divisor $A$. Let $X_i$ be the ambient variety of $\Aa_i$, $A_i$ the image of $A$ on $X_i$, and 
$$\lambda_i:=\inf\{s\geq 0\mid K_{\Aa_i}+sA_i\text{ is nef}/U\}$$
the scaling numbers. We may assume that this MMP does not terminate and let $\lim_{i\rightarrow+\infty}\lambda_i=\lambda$. There are two possibilities:

\medskip

\noindent\textbf{Case 1.} $\lambda<\lambda_i$ for any $i$. In this case, by \cite[Theorem 8.1]{Cas+25}, $\lambda=0$ and $K_{\Aa_i}+\lambda_iA_i$ is semi-ample$/U$ for any $i$. The theorem follows from Theorem \ref{thm: eomm implies tof with scaling} in this case.

\medskip

\noindent\textbf{Case 2.} $\lambda=\lambda_i>0$ for $i\gg 0$. In this case, $(\Aa_i,\lambda\overline{A})/U$ is a good minimal model of $(\Aa,\lambda\overline{A})/U$ for any $i\gg 0$. Let $h_i: X_i'\rightarrow X_i$ be a small $\mathbb Q$-factorialization of $X_i$ for each $i$ and let $\Aa_i':=h_i^*(\Aa_i,\lambda\overline{A})$ for any $i$. Then $\Aa_i'/U$ are $\mathbb Q$-factorial good minimal models of $(\Aa,\lambda\overline{A})/U$ for any $i$. By \cite[Theorem 2.5.2, Lemma 3.29]{Cas+25}, the induced birational map $\phi_{i,j}': X_i'\dashrightarrow X_j'$ is an isomorphism for some $0\ll i<j$. This is not possible because $\phi_{i,j}'$ is a non-trivial sequence of steps of a $K_{h_i^*\Aa_i}$-MMP.
\end{proof}

\begin{rem}
In Theorem \ref{thm: eomm implies tof with scaling}, by letting $t_0=t_1=0$, we indeed obtain the generalized pair analogue of \cite[Theorem 1.9]{Bir12}. Note that the NQC case was proven in \cite[Theorem 2.19]{TX24} but the non-NQC case remains open. 
\end{rem}

\begin{proof}[Proof of Theorem \ref{thm: scaling vs mmp main}]
    By \cite[Theorem 2.5.4]{Cas+25}, $(X,\Ff,t)$ has a KNQC minimal model. The theorem follows from Theorem \ref{thm: eo bswlcm implies eomm}.
\end{proof}

\section{Finiteness of ample models}\label{sec: Finiteness of ample models}

\subsection{Basic properties of ample models}

\begin{lem}
    Let $\pi: X\rightarrow U$ be a projective morphism between normal quasi-projective varieties. Let $D_1,D_2$ be two $\mathbb R$-Cartier $\mathbb R$-divisors on $X$.
    \begin{enumerate}
        \item For any birational map$/U$ $\phi: X\dashrightarrow X'$, if $\phi$ is a minimal model$/U$ (resp. weak lc model$/U$) of $D_1$ and $D_2$, then $\phi$ is also a minimal model$/U$ (resp. weak lc model$/U$) of $tD_1+(1-t)D_2$ for any $t\in [0,1]$.
        \item For any birational map$/U$ $\phi: X\dashrightarrow X'$, if $\phi$ is a minimal model$/U$ of $D_1$ and a weak lc model$/U$ of $D_2$, then $\phi$ is a minimal model of $tD_1+(1-t)D_2$ for any $t\in [0,1)$.
    \end{enumerate}
\end{lem}
\begin{proof}
Let $D_i':=\phi_*D_i$. Let $p: W\rightarrow X$ and $q: W\rightarrow X'$ be a resolution of indeterminacy of $\phi$ and write $p^*D_i=q^*D_i;+E_i$ for $i\in\{1,2\}$. Then $E_i\geq 0$, and $\Supp(p_*E_i)$ contains all prime divisors on $X$ that are exceptional$/X'$ if $\phi$ is a minimal model of $D_i$. We have
$$p^*(tD_1+(1-t)D_2)=q^*(tD_1'+(1-t)D_2')+\left(tE_1+(1-t)E_2\right)$$
and (1) and (2) follow.
\end{proof}

\begin{lem}\label{lem: semi-ample to ample}
    Let $\pi: X\rightarrow U$ be a projective morphism between normal quasi-projective varieties and let $D_1,D_2$ be two $\mathbb R$-divisors on $X$ such that $tD_1+(1-t)D_2$ is semi-ample$/U$ for any $t\in (0,1)$. Then:
    \begin{enumerate}
        \item There exists a contraction$/U$ $f: X\rightarrow Z$ such that $tD_1+(1-t)D_2=f^*H_t$ for some ample$/U$ $\mathbb R$-divisor $H_t$ for any $t\in (0,1)$.
        \item If $D_1$ is semi-ample$/U$, then there exists a contraction$/U$ $g: Z\rightarrow Y$ such that $D_1=(g\circ f)^*A$ for some ample$/U$ $\mathbb R$-divisor $A$.
    \end{enumerate}
\end{lem}
\begin{proof}
Let $f_t: X\rightarrow Z_t$ be the contraction induced by $tD_1+(1-t)D_2$, then the curves contracted by $f_t$ are exactly the curves $C$ such that $D_1\cdot C=0$ and $D_2\cdot C=0$, hence $f_t=f_{t'}$ for any $t,t'\in (0,1)$. (1) follows by letting $f:=f_t$. 

Moreover, when $D_1$ is semi-ample$/U$, since all curves $C$ contracted by $f_t$ satisfies that $D_1\cdot C=0$, we have $D_1=f^*D_Z$ for some $\mathbb R$-divisor $D_Z$ on $Z$.
When $D_1$ is semi-ample$/U$, $D_Z$ is semi-ample$/U$, so there exists a contraction$/U$ $g: Z\rightarrow Y$ such that $D_Z=g^*A$ for some ample$/U$ $\mathbb R$-divisor $A$. (2) follows.
\end{proof}

\begin{lem}\label{lem: semi-ample rationality}
     Let $\pi: X\rightarrow U$ be a projective morphism between normal quasi-projective varieties and let $D_1,D_2$ be two $\mathbb Q$-divisors on $X$ such that $tD_1+(1-t)D_2$ is semi-ample$/U$ for some irrational number $t\in (0,1)$. Let $f: X\rightarrow Z$ be the contraction induced by $tD_1+(1-t)D_2$. Then for any $s\in (t-\epsilon,t+\epsilon)$,
     $$sD_1+(1-s)D_2=f^*H_s$$
     for some ample$/U$ $\mathbb R$-divisor $H_s$. In particular, $sD_1+(1-s)D_2$ is semi-ample$/U$.
\end{lem}
\begin{proof}
By \cite[Lemma 5.3]{HLS24}, $sD_1+(1-s)D_2\sim_{\mathbb R,Z}0$ for any $s\in\mathbb R$. We have     
$$sD_1+(1-s)D_2=f^*H_s$$
for some $\mathbb R$-divisor $H_s$ for any $s\in\mathbb R$, and $H_t$ is ample$/U$. Thus $H_s$ is ample$/U$ for any $|s-t|\ll 1$.
\end{proof}

\subsection{Proof of the main theorems}

\begin{thm}\label{thm: technical finiteness of ample models}
Assume that $\Aa(t)/U:=(X,\Ff,B^{\ninv}+(1-t)B^{\inv},\Mm,t)/U$ is an algebraically integrable adjoint foliated structure for any $t\in [0,1]$. Let $0\leq s_0<s_1\leq 1$ be two real numbers such that 
\begin{itemize}
    \item $\Aa(t)/U$ is lc and has a bs-semi-ample model for any $t\in [s_0,s_1)$, and
    \item $\Aa(s_1)/U$ has a KNQC bs-weak lc model.
\end{itemize} 
Then there exist finitely many real numbers
$$s_0=:t_1<t_2<\dots<t_n:=s_1$$
satisfying the following. Let $\Ii:=\{\{t_i\},(t_i,t_{i+1})\mid 1\leq i\leq n-1\}$ if $\Aa(s_1)/U$ does not have bs-semi-ample model and let $\Ii:=\{\{t_i\},(t_j,t_{j+1})\mid 1\leq i\leq n, 1\leq j\leq n-1\}$ if $\Aa(s_1)/U$ has a bs-semi-ample model.

Then for any $\mathcal{P}\in\Ii$, there exists a birational map$/U$ $\phi_{\mathcal{P}}: X\dashrightarrow X_{\mathcal{P}}$ and a rational map$/U$ $\psi_{\mathcal{P}}: X\dashrightarrow Z_{\mathcal{P}}$ satisfying the following.
\begin{enumerate}
\item For any $\mathcal{P}\in\Ii$ and $t\in\mathcal{P}$, $\phi_{\mathcal{P}}$ is a $\mathbb Q$-factorial bs-good minimal model of $\Aa(t)/U$. 
\item For any $\mathcal{P}\in\Ii$ and $t\in\mathcal{P}$, $\psi_{\mathcal{P}}$ is the ample model$/U$ of $K_{\Aa(t)}$.
    \item For any $\mathcal{P},\mathcal{Q}\in\Ii$ such that $\mathcal{Q}\subset\partial\mathcal{P}$, there exist a unique contraction $\mu_{\mathcal{P},\mathcal{Q}}: Z_{\mathcal{P}}\rightarrow Z_{\mathcal{Q}}$ such that $\mu_{\mathcal{P},\mathcal{Q}}\circ\psi_{\mathcal{P}}=\psi_{\mathcal{Q}}$.
\end{enumerate}
Moreover:
\begin{enumerate}
    \item[(4)] If $X$ is $\mathbb Q$-factorial klt, then $\phi_{\mathcal{P}}$ does not extract any divisor. In particular, for any $\mathcal{P}\in\Ii$ and $t\in\mathcal{P}$, $\phi_{\mathcal{P}}$ is a $\mathbb Q$-factorial good minimal model of $\Aa(t)/U$. 
    \item[(5)] If $B$ is a $\mathbb Q$-divisor, $\Mm$ is a $\mathbb Q$-$\bb$-divisor, and $s_0,s_1$ are rational numbers, then we may choose $t_i$ so that they are rational numbers.
\end{enumerate}
\end{thm}

\begin{prop}\label{prop: ample model closed set case}
   Theorem \ref{thm: technical finiteness of ample models} holds if $\Aa(s_1)/U$ has a bs-semi-ample model.
\end{prop}
\begin{proof}
First we prove (1) and reduce (4) to (1).

Pick $t_0\in (s_0,s_1)$ and let $h: W\rightarrow X$ be a $\mathbb Q$-factorial qdlt modification of $\Aa(t_0)$. Then for any prime $h$-exceptional divisor $E$, $E$ is an lc place of $\Aa(t_0)$, hence an lc place of $\Aa(t)$ for any $t\in [0,1]$ as $\Aa(s_0)$ and $\Aa(s_1)$ are lc. We let $\Aa_W(t):=h^*\Aa(t)$ for any $t\in [0,1]$. By Lemma \ref{lem: mm preserved under dlt model}, possibly replacing $\Aa(t)$ with $\Aa_W(t)$ for any $t$, we may assume that $X$ is $\mathbb Q$-factorial klt. 

Assume that for any $t\in [s_0,s_1]$, there exists an open subset $V_t\ni t$ of $[s_0,s_1]$ and rational maps $\phi_{t}^-: X\dashrightarrow X_t^-$, $\phi_{t}^+: X\dashrightarrow X_t^+$, and $\phi_{t}: X\dashrightarrow X_t$ satisfying the following: For any $s\in V_t$, if $s<t$ (resp. $s=t$, $s>t$), then $\phi_t^-$ (resp. $\phi_t$, $\phi_t^+$) is a $\mathbb Q$-factorial good minimal model$/U$ of $\Aa(s)/U$. Then $\{V_t\}_{t\in [s_0,s_1]}$ is an open covering of $[s_0,s_1]$, hence there exist finitely many real numbers
$$s_0=:t_1<t_2<\dots<t_n:=s_1,$$
such that $\{V_{t_i}\}_{1\leq i\leq n}$ is an open covering of $[s_0,s_1]$, and that $V_{t_i}\not\subset V_{t_j}$ for an $i\not=j$. In particular, $V_{t_i}\cap V_{t_{i+1}}\not=\emptyset$,
hence $\phi_{t_i}^+=\phi_{t_{i+1}}^-$ for any $i$. These $t_i$ satisfy our requirements. Therefore, we only need to find $V_t$ and $\phi_{t}^-,\phi_t,\phi_t^+$ as above for any $t\in [s_0,s_1]$.

Fix $t\in [s_0,s_1]$. By Theorem \ref{thm: eo bswlcm implies eomm}, we may run a $K_{\Aa(t)}$-MMP$/U$ with scaling of an ample divisor which terminates with a $\mathbb Q$-factorial good minimal model $\Aa'(t)/U$ of $\Aa(t)/U$. Let $\phi_t: X\dashrightarrow X'=:X_t$ be the induced birational map and let $\Aa'(s):=(\phi_t)_*\Aa(s)$ for any $s\in [0,1]$. Then $X'$ is klt, and there exists a contraction$/U$ $f: X'\rightarrow Z$ and an ample$/U$ $\mathbb R$-divisor $H$ on $Z$ such that
$$K_{\Aa'(t)}=f^*H.$$
There exists $\epsilon_1>0$, such that for any $s\in I_1:=[t-\epsilon_1,t+\epsilon_1]\cap [s_0,s_1]$, $\phi_t$ is also a sequence of steps of a $K_{\Aa(s)}$-MMP$/U$, hence $\phi_t$ is $K_{\Aa(s)}$-negative and $K_{\Aa'(s)}$ is lc. By Lemma \ref{lem: minimal model same after running mmp}, $\Aa'(s)/U$ has a good minimal model for any $s\in I_1$. By Lemma \ref{lem: trivial mmp}, there exists $\epsilon_2\in (0,\epsilon_1)$, such that for any $s\in I_2:=[t-\epsilon_2,t+\epsilon_2]\cap [s_0,s_1]$, any sequence of steps of a $K_{\Aa'(s)}$-MMP$/U$ is $K_{\Aa'(t)}$-trivial, 
hence a sequence of steps of a $K_{\Aa'(s)}$-MMP$/Z$.

Let $t_-:=\max\{t-\epsilon_2,s_0\}$ and $t_+:=\min\{t+\epsilon_2,s_1\}$. By Theorem \ref{thm: eo bswlcm implies eomm}, we may run a $K_{\Aa'(t_-)}$-MMP$/U$ (resp. $K_{\Aa'(t_+)}$-MMP$/U$) which terminates with a good minimal model $\Aa_-(t_-)/U$ (resp. $\Aa_+(t_+)/U$) of $\Aa'(t_-)/U$ (resp.$\Aa'(t_+)/U$) with induced birational map $\phi_t^-: X'\dashrightarrow X_t^-$ (resp. $\phi_t^+: X'\dashrightarrow X_t^+$).  $\phi_t^-,\phi_t^+$ are maps$/Z$. Let $\Aa_-(s):=(\phi_t^-)_*\Aa'(s)$ and $\Aa_+(s):=(\phi_t^+)_*\Aa'(s)$ for any $s\in [0,1]$. Then for $s\in [t_-,t)$ (resp. $(t,t_+]$), $K_{\Aa'(s)}\sim_{\mathbb R,Z} \mu_s K_{\Aa(t_-)}$ (resp. $\sim_{\mathbb R,Z} \mu_s K_{\Aa(t_+)}$) for some $\mu_s>0$. Thus $\phi_t^-$ (resp. $\phi_t^+$) is also a sequence of steps of a $K_{\Aa(s)}$-MMP$/Z$ and $K_{\Aa_-(s)}$ (resp. $K_{\Aa_+(s)}$) is nef$/Z$. Since any sequence of steps of a $K_{\Aa'(s)}$-MMP$/U$ is a sequence of steps of a $K_{\Aa'(s)}$-MMP$/Z$, $K_{\Aa_-(s)}$ (resp. $K_{\Aa_+(s)}$) is nef$/U$, so $\Aa_-(s')/U$ (resp. $\Aa_+(s')/U$)  is a minimal model of $\Aa'(s)/U$. By Lemma \ref{lem: minimal model same after running mmp}, $\Aa_-(s)/U$ (resp.  $\Aa_+(s)/U$) is a minimal model of $\Aa(s)/U$. By Lemma \ref{lem: g-pair version bir12 2.7}, $K_{\Aa_-(s)}$ (resp.  $K_{\Aa_+(s)}$ ) is semi-ample$/U$. We may let $V_t:=(t_-,t_+)$ if $t\not\in\{s_0,s_1\}$, $V_{s_0}=[s_0,s_{0,+})$, and $V_{s_1}:=(s_{1,-},s_1]$. (1) follows by our construction of $\phi_t^-,\phi_t,\phi_t^+$.

(2) follows from (1), the uniqueness of ample model, and Lemma \ref{lem: semi-ample to ample}(1). (3) follows from Lemma \ref{lem: semi-ample to ample}(2).

Finally, we prove (5). Assume that $t_i\not\in\mathbb Q$ for some $i$, then $i\not\in\{1,n\}$. Let $\phi_+:=\phi_{(t_i,t_{i+1})}$, $\phi_-:=\phi_{(t_{i-1},t_{i})}$, $X^+:=X_{(t_i,t_{i+1})}$, and $X^-:=X_{(t_{i-1},t_{i})}$.

By comparing discrepancies and since nef is a closed condition,
$\phi_+$ and $\phi_-$ are $\mathbb Q$-factorial bs-weak lc models of $\Aa(t_i)/U$. For any $t\in [0,1]$, write
$$\Aa_+(t):=\left((\phi_+)_*\Aa(s),\sum_D(t\epsilon_{\Ff}(D)+(1-t))D\right)$$
where the sum runs through all prime divisors $D$ on $X^+$ that are exceptional$/X$ and 
$$\Aa_-(t):=\left((\phi_-)_*\Aa(t),\sum_D(t\epsilon_{\Ff}(D)+(1-t))D\right)$$
where the sum runs through all prime divisors $D$ on $X^-$ that are exceptional$/X$. 

For any prime divisor $D$ on $X$ that is exceptional$/X^+$ (resp. $X^-$), we have that $a(D,\Aa(t_i))\leq a(D,\Aa_+(t_i))$ (resp. $a(D,\Aa(t_i))\leq a(D,\Aa_-(t_i))$). Assume that $a(D,\Aa(t_i))=a(D,\Aa_+(t_i))$ (resp. $a(D,\Aa(t_i))=a(D,\Aa_-(t_i))$). Since $t_i$ is irrational, we have that $a(D,\Aa(t))=a(D,\Aa_+(t))$ (resp. $a(D,\Aa(t))=a(D,\Aa_-(t))$) for any $t\in [0,1]$, which is not possible because $a(D,\Aa(t))<a(D,\Aa_+(t))$ (resp.  $a(D,\Aa(t))<a(D,\Aa_-(t))$ when $t\in (t_i,t_{i+1})$ (resp. $t\in (t_{i-1},t_i)$). Therefore, $a(D,\Aa(t_i))<a(D,\Aa_+(t_i))$ (resp. $a(D,\Aa(t_i))<a(D,\Aa_-(t_i))$). Thus $\phi_+$ and $\phi_-$ are both bs-good minimal models of $\Aa(t_i)$.

By Lemma \ref{lem: semi-ample rationality}, $K_{\Aa_+(t)}$ and $K_{\Aa_-(t)}$ are semi-ample$/U$ for any $|t-t_i|\ll 1$. Since there are only finitely many prime divisors $D$ on $X$ that are exceptional$/X^+$ or $X^-$, there exist rational numbers $t_i^-,t_i^+$ such that $t_i^-<t_i<t_i^+$, and $\Aa_+(t)/U$ and $\Aa_-(t)/U$ are $\mathbb Q$-factorial good minimal models of $\Aa(t_i)/U$ for any $t\in [t_i^-,t_i^+]$. By Lemma \ref{lem: semi-ample rationality} again, for any $|t-t_i|\ll 1$, $\Aa_+(t)/U$ and $\Aa_-(t)/U$ have the same ample model. By the uniqueness of ample models, $\psi_{\{t_i\}}=\psi_{(t_{i-1},t_i)}=\psi_{(t_{i},t_{i+1})}$. Now we may replace $\Ii$ by
$$\{\{t_j\},\{t_i^-\},\{t_i^+\},(t_k,t_{k+1}),(t_{i-1},t_i^-),(t_i^-,t_i^+),(t_i^+,t_{i+1})\mid j\not=i,k\not=i,i-1\}.$$
let $$\psi_{(t_{i-1},t_i^-)}:=\psi_{\{t_i^-\}}:=\psi_{(t_{i}^-,t_i^+)}:=\psi_{\{t_{i}^+\}}:=\psi_{(t_{i}^+,t_{i+1})}:=\psi_{t_i},$$
$$\phi_{(t_{i-1},t_i^-)}:=\phi_{\{t_i^-\}}:=\phi_{(t_{i}^-,t_i^+)}:=\phi_{\{t_{i}^+\}}:=\phi_{(t_{i-1},t_i)},$$
and
$$\phi_{(t_{i}^+,t_{i+1})}:=\phi_{(t_{i},t_{i+1})}.$$
We are done by induction on the number of irrational numbers which appear in $\{t_i\}_{i=1}^n$.
\end{proof}

\begin{lem}\label{lem: when t goes to 1}
    Let $\Aa(t)/U:=(X,\Ff,B^{\ninv}+(1-t)B^{\inv},\Mm,t)/U$ be an algebraically integrable adjoint foliated structure for any $t\in [0,1]$. Assume that
    \begin{enumerate}
        \item $\Aa(t_0)/U$ is lc and has a KNQC bs-weak lc model for some $t_0\in (0,1]$,
        \item either $t_0<1$, or $\Aa(t)$ is lc for some $t\in [0,1)$, and
        \item there exists a strictly increasing sequence $\{t_i\}_{i=1}^{+\infty}$ such that $\lim_{i\rightarrow+\infty}t_i=t_0$ and $\Aa(t_i)/U$ has a KNQC bs-weak lc model.
    \end{enumerate}
    Then there exists $\epsilon_0>0$ and a birational map $\phi: X\dashrightarrow X'$ satisfying the following: For any $t\in (t_0-\epsilon_0,t_0)$, $\phi$ is a $\mathbb Q$-factorial KNQC bs-minimal model of $\Aa(t)/U$. Moreover, if $X$ is $\mathbb Q$-factorial klt, then $\phi$ does not extract any divisor.
\end{lem}
\begin{proof}
By \cite[Proposition 3.3]{Cas+24}, $\Aa(t)$ is lc for any $t\in [0,t_0]$. Let $h: W\rightarrow X$ be a $\mathbb Q$-factorial qdlt modification of $\Aa(t)$ for some $t\in (0,t_0)$. Then any prime $h$-exceptional divisor $D$ is an lc place of $\Aa(t)$, hence an lc place of $\Aa(s)$ for any $s\in [0,1]$ as $\Aa(0)$ and $\Aa(t_0)$ are lc. We let $\Aa_W(t):=h^*\Aa(t)$ for any $t\in [0,1]$. By Lemmas \ref{lem: same weak glc model under pullback} and \ref{lem: mm preserved under dlt model}, we may replace $\Aa(t)$ with $\Aa_W(t)$ and assume that $X$ is $\mathbb Q$-factorial klt.

By Theorem \ref{thm: eo bswlcm implies eomm}, we may run a $K_{\Aa(t_0)}$-MMP$/U$ which terminates with a $\mathbb Q$-factorial KNQC minimal model $\Aa'(t_0)/U$ of $\Aa(t_0)/U$ with induced birational map $\phi: X\dashrightarrow X'$. Let $\Aa'(t):=\phi_*\Aa(t)$ for any $t\in [0,1]$. Then there exists $t_0>\epsilon_1>0$, such that $\phi$ is also a sequence of steps of a $K_{\Aa(t_0-\epsilon)}$-MMP$/U$ for any $\epsilon\in [0,\epsilon_1]$. By Lemma \ref{lem: minimal model same after running mmp}, any minimal model of $\Aa'(t_0-\epsilon)/U$ is a minimal model of $\Aa(t_0-\epsilon)/U$. By Lemma \ref{lem: trivial mmp}, there exists $0<\epsilon_2<\epsilon_1$, such that for any $\epsilon\in [0,\epsilon_2]$, any sequence of steps of a $K_{\Aa'(t_0-\epsilon)}$-MMP$/U$ is $K_{\Aa'(t_0)}$-trivial. 

Pick $i\gg 0$ such that $t_i\in (t_0-\epsilon_2,t_0)$. By Theorem \ref{thm: eo bswlcm implies eomm}, we may run a $K_{\Aa'(t_i)}$-MMP$/U$ which terminates with a $\mathbb Q$-factorial KNQC minimal model $\Aa''(t_i)/U$ of $\Aa'(t_i)/U$ with induced birational map $\psi: X'\dashrightarrow X''$. Let $\Aa''(t):=\psi_*\Aa'(t)$ for any $t\in [0,1]$. Since $\psi$ is $K_{\Aa'(t_0)}$-trivial, $\psi$ is also a sequence of steps of $K_{\Aa'(t)}$-MMP$/U$ for any $t\in [t_i,t_0)$. Moreover, since $K_{\Aa''(t_0)}$ and $K_{\Aa''(t_i)}$ are both NQC$/U$, $K_{\Aa''(t)}$ is NQC$/U$ for any $t\in [t_i,t_0)$. Thus $\Aa''(t)/U$ is a $\mathbb Q$-factorial KNQC minimal model of $\Aa'(t)/U$, hence a $\mathbb Q$-factorial KNQC minimal model of $\Aa(t)/U$. We may take $\epsilon_0:=t_0-t_i$.
\end{proof}

\begin{proof}[Proof of Theorem \ref{thm: technical finiteness of ample models}]
By Proposition \ref{prop: ample model closed set case}, we may assume that $\Aa(s_1)/U$ does not have a bs-semi-ample model.

 (1)  By Lemmas \ref{lem: when t goes to 1} and \ref{lem: g-pair version bir12 2.7}, there exists $s_1-s_0>\epsilon>0$ and a birational map $\phi: X\dashrightarrow X'$, such that $s_1-\epsilon$ is rational, and for any $t\in (s_1-\epsilon,s_1)$, $\phi$ is a $\mathbb Q$-factorial bs-good minimal model of $\Aa(t)/U$. By Lemma \ref{lem: semi-ample to ample}(1), there exists a contraction $f: X'\rightarrow Z$ such that the induced birational map $\psi: f\circ\phi$ is the ample model$/U$ of $K_{\Aa(t)}$.
 
 Since $\Aa(s_1-\epsilon_0)/U$ has a bs-semi-ample model, by Proposition \ref{prop: ample model closed set case}, there exist finitely many real numbers
 $$s_0=:t_1<t_2<\dots<t_{n-1}:=s_1-\epsilon$$
satisfying the following. Let $\Ii':=\{\{t_i\},(t_j,t_{j+1})\mid 1\leq i\leq n, 1\leq j\leq n-1\}$. Then for any $\mathcal{P}\in\Ii'$, there exists a birational map$/U$ $\phi_{\mathcal{P}}: X\dashrightarrow X_{\mathcal{P}}$, and a rational map$/U$ $\psi_{\mathcal{P}}: X\dashrightarrow Z_{\mathcal{P}}$ satisfying the following.
\begin{itemize}
\item For any $t\in\mathcal{P}$, $\phi_{\mathcal{P}}$ is a $\mathbb Q$-factorial bs-good minimal model of $\Aa(t)/U$. 
\item For any $\mathcal{P}\in\Ii'$ and $t\in\mathcal{P}$, $\psi_{\mathcal{P}}$ is the ample model of $\Aa(t)/U$.
    \item For any $\mathcal{P},\mathcal{Q}\in\Ii'$ such that $\mathcal{Q}\subset\partial\mathcal{P}$, there exist a unique contraction $\mu_{\mathcal{P},\mathcal{Q}}: Z_{\mathcal{P}}\rightarrow Z_{\mathcal{Q}}$ such that $\mu_{\mathcal{P},\mathcal{Q}}\circ\psi_{\mathcal{P}}=\psi_{\mathcal{Q}}$.
    \item If $X$ is $\mathbb Q$-factorial klt, then $\phi_{\mathcal{P}}$ does not extract any divisor. In particular, for any $\mathcal{P}\in\Ii'$ and $t\in\mathcal{P}$, $\phi_{\mathcal{P}}$ is a $\mathbb Q$-factorial good minimal model of $\Aa(t)/U$. 
    \item If $B$ is a $\mathbb Q$-divisor, $\Mm$ is a $\mathbb Q$-$\bb$-divisor, and $s_0$ is a rational number, then we may choose $t_i$ so that they are rational numbers.
\end{itemize}
We may let $t_n:=s_1$, $\Ii:=\Ii'\cup\{(s_1-\epsilon,s_1)\}$, $\phi_{(s_1-\epsilon,s_1)}:=\phi$, and $\psi_{(s_1-\epsilon,s_1)}:=\psi$. Then (1)(2)(4)(5) follow. By comparing discrepancies and since nef$/U$ is a closed condition, $\phi$ is a $\mathbb Q$-factorial bs-semi-ample model of $\Aa(s_1-\epsilon)/U$. Thus (3) follows from Lemma \ref{lem: semi-ample to ample}.
\end{proof}

\begin{proof}[Proof of Theorem \ref{thm: main canonical F finiteness of ample models f general type}]
Since $X$ is klt, possibly replacing $X$ with a small $\mathbb Q$-factorialization, we may assume that $X$ is $\mathbb Q$-factorial. Since $\Aa_{\lambda}$ is of general type, by \cite[Theorem 5.1]{Cas+25}, $K_{\Ff}$ is pseudo-effective. Let $h: (X',\Ff',B')\rightarrow (X,\Ff,0)$ be a $\mathbb Q$-factorial ACSS modification of $\Ff$ such that $(X',\Ff',B';G)/Z$ is ACSS associated with contraction $f: X'\rightarrow Z$. Since $\Ff$ is canonical, $h$ is an isomorphism over the generic point of $Z$ and $B'=0$. Let $F$ be a general fiber of $f$, then $F$ is klt and $K_F=h^*K_{\Ff}|_F=h^*K_X|_F$ is big, so $F$ has a good minimal model. By \cite[Theorem 3.15]{HJLL24}, $\Ff'$ has a minimal model$/Z$. By the cone theorem \cite[Theorem 3.9]{ACSS21}, $\Ff'$ has a minimal model. By Lemma \ref{lem: minimal model same after running mmp}, $\Ff$ has a minimal model, which is KNQC as $K_{\Ff}$ is a $\mathbb Q$-divisor. 

Pick $0<\delta\ll 1$ such that $\lambda-\delta\in\mathbb Q$ and let $\epsilon:=\max\{0,\lambda-\delta\}$. Then $\Aa_t$ is klt and of general type for any $t\in [\epsilon,1)$, hence $(X,\Ff,t)$ has a good minimal model for any $t\in [\epsilon,1)$ by \cite[Theorem 2.1.1]{Cas+25}. The theorem follows from Theorem \ref{thm: technical finiteness of ample models}.
\end{proof}

\end{document}